\documentclass[10pt, a4paper]{article}
\date{4.2.2003}
\author{Frank Schuhmacher\footnote{Supported by:
Doktorandenstipendium vom Deutschen Akademischen Austauschdienst im Rahmen
des gemeinsamen Hochschulsonderprogramms III des Bundes und der
L\"ander}}

\title{Hochschild Cohomology for Complex Spaces and Noetherian Schemes}

\usepackage{a4wide}
\usepackage[intlimits]{amsmath}
\usepackage{amssymb}
\usepackage{theorem} 
\usepackage{epic}
\usepackage{amscd}
\usepackage[all]{xypic}
\usepackage{makeidx}

\newcounter{punkt}
\newenvironment{liste}{\begin{list}{(\roman{punkt})}
{\usecounter{punkt} 
\setlength{\topsep}{0cm}
\setlength{\itemsep}{0.1cm}
\setlength{\parsep}{0cm}
}}{\end{list}}

\newcounter{defi}[section]
\newenvironment{defi}{\refstepcounter{defi}\vspace{3pt}\hfill\\
 \textbf{Definition }\textbf{\arabic{section}.\arabic{defi}}
 }
 {\vspace{1pt}\hfill\\}

\newcounter{satz}[section]
\newenvironment{satz}{\refstepcounter{satz}\vspace{3pt}\hfill\\
 \textbf{Theorem }\textbf{\arabic{section}.\arabic{satz}}
 \begin{it}}
 {\end{it}\vspace{1pt}\hfill\\}

\newcounter{prop}[section]
\newenvironment{prop}{\refstepcounter{prop}\vspace{3pt}\hfill\\
 \textbf{Proposition }\textbf{\arabic{section}.\arabic{prop}}
 }
 {\vspace{1pt}\hfill\\}

\newcounter{lemma}[section]
\newenvironment{lemma}{\refstepcounter{lemma}\vspace{3pt}\hfill\\
 \textbf{Lemme }\textbf{\arabic{section}.\arabic{lemma}}
 \begin{it}}
 {\end{it}\vspace{1pt}\hfill\\}

\newcounter{kor}[section]
\newenvironment{kor}{\refstepcounter{kor}\vspace{3pt}\hfill\\
 \textbf{Corollary }\textbf{\arabic{section}.\arabic{kor}}
 \begin{it}}
 {\end{it}\vspace{1pt}\hfill\\}

\newcounter{bem}[section]
\newenvironment{bem}{\refstepcounter{bem}\vspace{3pt}\hfill\\
 \textbf{Remark }\textbf{\arabic{section}.\arabic{bem}}
 }
 {\vspace{1pt}\hfill\\}

\newcounter{beisp}[section]
\newenvironment{beisp}{\refstepcounter{beisp}\vspace{3pt}\hfill\\
 \textbf{Example }\textbf{\arabic{section}.\arabic{beisp}}
 }
 {\vspace{1pt}\hfill\\}

\newenvironment{axiom}{\vspace{3pt}\hfill\\
 \textbf{Axiom }
 }
 {\vspace{1pt}\hfill\\}

\newenvironment{axioms}{\vspace{3pt}\hfill\\
 \textbf{Axioms }}
 {\vspace{1pt}\hfill\\}
\renewcommand{\labelenumi}{(\roman{enumi})}

\newcommand{\impl}{\Rightarrow }
\newcommand{\incl}{\hookrightarrow}

\newcommand{\nach}{\longrightarrow}
\newcommand{\sub}{\subseteq}

\newcommand{\isom}{\cong}
\newcommand{\lan}{\langle}
\newcommand{\ran}{\rangle}

\newcommand{\ZZ}{\mathbb{Z}}

\newcommand{\QQ}{\mathbb{Q}}
\newcommand{\CC}{\mathbb{C}}
\newcommand{\LL}{\mathbb{L}}
\newcommand{\KK}{\mathbb{K}}

\newcommand{\sF}{\mathcal{F}}
\newcommand{\sG}{\mathcal{G}}

\newcommand{\sI}{\mathcal{I}}

\newcommand{\Oh}{\mathcal{O}}

\newcommand{\sR}{\mathcal{R}}

\newcommand{\ppp}{\cdot\ldots\cdot}

\newcommand{\Mod}{\mathfrak{Mod}}
\newcommand{\Alg}{\mathfrak{Alg}}

\newcommand{\sN}{\mathcal{N}}
\newcommand{\sA}{\mathcal{A}}
\newcommand{\sS}{\mathcal{S}}
\newcommand{\sB}{\mathcal{B}}

\newcommand{\Hk}{\mathbb{H}}
\newcommand{\cone}{\operatorname{cone}}
\newcommand{\dach}{\wedge}
\newcommand{\Ext}{\operatorname{Ext}}
\newcommand{\sT}{\mathcal{T}}
\newcommand{\HOM}{\mathcal{H}\!\mathit{om}}
\newcommand{\EXT}{\mathcal{E}\!\mathit{xt}}
\newcommand{\MOD}{\mathcal{M}\!\mathit{od}}
\newcommand{\COH}{\mathcal{C}\!\mathit{oh}}
\newcommand{\grm}{\gr(\sM)}
\newcommand{\grc}{\gr(\sC)}
\newcommand{\grmn}{\gr(\sM)^{\sN}}
\newcommand{\grcn}{\gr(\sC)^{\sN}}
\newcommand{\cycl}{\operatorname{cycl}}
\newcommand{\sX}{\mathcal{X}}
\newcommand{\im}{\operatorname{im}}
\renewcommand{\Im}{\operatorname{Im}}
\newcommand{\sM}{\mathcal{M}}
\newcommand{\sC}{\mathcal{C}}
\newcommand{\X}{\mathcal{X}}
\newcommand{\Y}{\mathcal{Y}}

\newcommand{\HH}{\operatorname{HH}}

\newcommand{\N}{\operatorname{N}}
\renewcommand{\P}{\operatorname{P}}

\newcommand{\T}{\mathsf{T}}

\newcommand{\m}{\mathfrak{m}}

\renewcommand{\HH}{\operatorname{HH}}
\newcommand{\bew}{\bf Proof: \rm}
\newcommand{\qed}{\hfill $ \square$ \\}
\newcommand{\spec}{\operatorname{Spec}}

\newcommand{\Hom}{\operatorname{Hom}}
\newcommand{\Mult}{\operatorname{Mult}}
\newcommand{\Tor}{\operatorname{Tor}}

\newcommand{\id}{\operatorname{Id}}

\newcommand{\gr}{\operatorname{gr}}
\newcommand{\tot}{\operatorname{tot}}
\newcommand{\ba}{\operatorname{bar} }
\newcommand{\ob}{\operatorname{ob} }
\newcommand{\Kern}{\operatorname{Kern} }
\renewcommand{\ker}{\operatorname{kern} }
\newcommand{\Kokern}{\operatorname{Cokern} }
\newcommand{\kokern}{\operatorname{cokern} }
\newcommand{\lz}{\hfill\newline}

\newcommand{\ot}{\otimes}

\begin{document}

\maketitle

\begin{abstract}

The classical HKR-theorem gives an isomorphism of the n-th Hochschild
cohomology of a smooth algebra and the n-th exterior power of its
module of K\"ahler differentials. Here we generalize it for simplicial,
graded and anticommutative objects in ``good pairs of categories''.
We apply this generalization to complex spaces and noetherian schemes
and deduce two decomposition theorems
for  their (relative) Hochschild cohomology (special cases of those were
recently shown by Buchweitz-Flenner and Yekutieli). The first one
shows that Hochschild cohomology contains tangent cohomology:
$\HH^n(X/Y,\sM)=\coprod_{i-j=n}\Ext^i(\dach^j\LL(X/Y),\sM)$.
The left side is the n-th Hochschild cohomology of $X$ over $Y$
with values in $\sM$. The right hand-side
contains the $n$-th relative
tangent cohomology $\Ext^n(\LL(X/Y),\sM)$
as direct factor. 
The second consequence is a decomposition theorem
for Hochschild cohomology of complex analytic manifolds
and smooth schemes in characteristic zero:
$\HH^n(X)=\coprod_{i-j=n}H^i(X,\dach^j\sT_X).$
On the right hand-side we have the sheaf cohomology of the exterior
powers of the tangent complex.\\


\end{abstract}

\tableofcontents

\section*{Introduction}
A better title for this paper would be:
``Hochschild Cohomology for Admissible Pairs of
Categories and Application to Complex Spaces and Noetherian
Schemes''. Since this title would be too lengthy, 
and admissible pairs of categories
seem not to be so well known, I did not mention them in the title. 
Admissible pairs of categories are pairs $(\sC,\sM)$, where
$\sC$ is a certain category of algebras (one should think of global sections
of the structure sheaf of an ``affine'' space\footnote{In the analytic
  context, by an ``affine'' space we mean a Stein compact.})
and $\sM$ is a category of certain modules over objects in $\sC$
(one should think of global sections of coherent modules over an affine). 
They were invented by Bingener and
Kosarew in \cite{BinKos} and are quite useful in deformation theory,
since the PO-algebras and modules of Palamodov also fit into
this definition.  
To describe spaces globally, one has to consider simplicial objects
in $\sC$ and $\sM$ for an admissible pair $(\sC,\sM)$, 
i.e. functors from the ``nerf'' of an affine
covering of the space to $\sC$ and $\sM$. 
Each construction, using
admissible pairs, which is canonical in a sense, can be generalized to
the simplicial case. To generalize non-canonical constructions one has
to do some work.

In this paper,we use admissible pairs of categories to unify
the algebraic Hochschild theory, that can be found in several 
textbooks (for ex.\cite{Loday}), and the geometrical approach, 
which in the analytical
context is due to Buchweitz and Flenner.
The central result will be a generalization of the
classical Hochschild-Kostant-Rosenberg theorem for
differential graded algebras. From this generalization we
will deduce the following HKR-type theorem:\\
When $X\nach Y$ is a morphism of complex spaces
(paracompact and separated)
or a separated morphism of finite type of noetherian schemes 
in characteristic zero,
then there is a quasi-isomorphism
\begin{equation}\label{HKRT}
\Hk(X/Y)\approx\dach\LL(X/Y)
\end{equation} 
over $\Oh_X$, where $\Hk(X/Y)$
is the relative Hochschild complex of $X$ over $Y$ 
(see section~\ref{Appl})
and $\LL(X/Y)$ is the relative cotangent complex.
This statement is also true, when $Y$ is just a single point and
$X$ a smooth scheme in characteristic zero.
From this main result we will deduce two nice decomposition theorems
for  Hochschild cohomology. The first one
shows that Hochschild cohomology contains tangent cohomology:
\begin{equation}\label{HUNT}
\HH^n(X/Y,\sM)=\coprod_{i-j=n}\Ext^i(\dach^j\LL(X/Y),\sM).
\end{equation}
The left side is the n-th Hochschild cohomology of $X$ over $Y$
with values in $\sM$.
The right hand-side
contains the $n$-th relative
tangent cohomology $\Ext^n(\LL(X/Y),\sM)$
as direct factor. For complex spaces, this decomposition,
as well as equation~\ref{HKRT}, was already shown in
a different way by Buchweitz and Flenner \cite{BuFl2}.\\

The second consequence is a decomposition theorem
for Hochschild cohomology of complex analytic manifolds
and smooth schemes in characteristic zero:
$$\HH^n(X)=\coprod_{i-j=n}H^i(X,\dach^j\sT_X).$$
On the right hand-side we have the sheaf cohomology of the exterior
powers of the tangent complex.
A proof of this result for complex analytic manifolds
was announced (but not yet published) by Kontsevich.
For schemes, the latter result is also 
proven by Yekutieli \cite{Yeku}.

To make this paper more
independent, we will state all axioms and definitions that
we need about admissible pairs of categories in section
~\ref{Apoc}. One reason for this is, that the main reference
\cite{BinKos} is not any more available.

To prove a global statement as equation~\ref{HKRT}, we need 
three steps: First, we have to prove it for good pairs
of categories. Secondly, we have to generalize it to simplicial
objects in $\sC$ and $\sM$. Third, by using a Cech-construction,
we prove them for sheaves of algebras or modules.\\

This paper is organized as follows:
In the beginning of section~\ref{Apoc}, we state the definitions 
for admissible pairs of categories. In \ref{fmaa}, we will precise the notation
of free objects and of good pairs of categories. In \ref{ago} and
\ref{dgo}, we explain how to construct, starting from an admissible
(resp. good) pair $(\sC,\sM)$ of categories, the admissible (resp. good)
pairs $(\grc,\grm)$ and $(\gr^2(\sC),\gr^2(\sM))$, containing 
(anti-commutative) graded and double graded objects. In \ref{reso}
we define resolutions and resolvents and prove a first important
statement, saying that two resolvents of a homomorphism
of (simplicial) algebras are homotopy equivalent. This is an
improvement
of \cite{BinKos}, proposition (8.4), saying that two resolvents are
quasi-isomorphic.
In \ref{barr}, we generalize the definition of the (cyclic) bar
complex for admissible pairs of categories and state their main
properties. (In the classical literature, the cyclic bar complex is 
called Hochschild complex.)
In \ref{hans}, we adapt the definition of regular sequences
to the graded anticommutative context. Here we prove the equivalence
of four different conditions.
In \ref{umod}, we have to prove several statements, 
concerning the universal module of differentials.
The reader who is only interested in the theory of schemes or
complex spaces, i.e. in the examples $(\sC^{(0)},\sM^{(0)})$
and $(\sC^{(1)},\sM^{(1)})$ of good pairs (see example~\ref{Apoc}.\ref{exop}),
can leave out the lecture of this subsection, since for this examples,
all statements proven here, are well-known.
In section~\ref{hchc}, we define Hochschild complexes and
Hochschild cohomology for admissible pairs of categories.
In~\ref{tac} we show that the definitions in \ref{secta} generalize
the classical definitions in the algebraic context. 
In section~\ref{Adt} we prove the main result, i.e. a generalization
of the classical HKR-theorem and deduce a decomposition theorem
for Hochschild cohomology.
In section~\ref{Appl}, this results are applied to schemes and
complex spaces. For this, we have to introduce simplicial technics and to
define (in \ref{hacop}) resolvents of homomorphisms of complex spaces and
schemes. Then we can define Hochschild complexes and Hochschild
cohomology in the same way as Buchweitz and Flenner \cite{BuFl2}, 
and the announced decomposition theorems will
be deduced by the generalized HKR-theorem in the last two subsections.\\

\textbf{Acknowledgments:} I want to express my gratitude to Siegmund
Kosarew for many useful suggestions and to Ragnaz-Olaf Buchweitz and
Hubert Flenner for a preliminary version of their paper on Hochschild 
cohomology.

\textbf{Conventions:}
For a ring $k$, we denote the category of $k$-modules by $k$-$\Mod$. 
When we work with a morphism $f:A\nach B$ in any category, we denote
the kernel of $f$ in the cetegorial sense by $\ker f$.
So $\ker f$ is a morphism $K\nach A$, where $K$ is an object, determined up to
a canonical isomorphism. By $\Kern f$, we mean the object $K$.
In the same manner we use the notations $\kokern$, $\Kokern$, $\im$ and $\Im$.
So, for example, we have $\Im f=\Kern(\kokern f)$. 
We will use $\approx$, to denote quasi-isomorphisms and 
$\simeq$ to denote homotopy-equivalences. 
We will use the letter $D$ to denote derived categories and 
$K$ to denote homotopy categories, i.e. the localization of categories
by homotopy-equivalences. 
\pagebreak

\section{Admissible pairs of categories}\label{Apoc}

We fix a ground ring $\KK$ (in our main reference \cite{BinKos}
$\KK$ is the field $\QQ$, so here we start with a more general setting).
Denote by $\sC$ a category of commutative $\KK$-algebras
and by $\sC$-$\Mod$ the category of all modules over algebras in $\sC$
and let $\sM$ be a subcategory of $\sC$-$\Mod$. Then the pair
$(\sC,\sM)$ is called an \textbf{admissible pair of categories}
if the following conditions hold:
\renewcommand{\labelenumi}{(\arabic{enumi})}
\renewcommand{\labelenumii}{(\arabic{enumi}.\arabic{enumii})}
\begin{enumerate}
\item
In $\sC$ there exist finite fibered sums, that we denote as usual 
by $A\ot^{\sC}_kB$.
\item
When $\phi:A\nach B$ is a homomorphism in $\sC$ and $N$ a module
in $\sM(B)$, then $N$ is via $\phi$ an object of $\sM(A)$
and for each module $M$ in $\sM(A)$, $\Hom_{\sM(A)}(M,N)$
is the set of all $\phi$-homomorphisms $M\nach N$ in $\sM$.
\item
Let $A$ be an algebra in $\sC$. Then $\sM(A)$ is an additive category,
in which kernels and cokernels exist. Further
$\sC_A$ is a subcategory of $\sM(A)$ and the functor of $\sM(A)$
in $A$-$\Mod$ commutes with kernels and finite direct sums.
\item
Let $\phi:A\nach B$ a homomorphism in $\sC$ and $u:M\nach N$
a homomorphism in $\sM(B)$. Let $L$ resp. $L'$ be the kernel
of $u$ resp. $u_{[\phi]}$ in $\sM(B)$ resp. $\sM(A)$.
Then the canonical map $L'\nach L_{[\phi]}$ is an isomorphism
in $\sM(A)$.
\item\label{zp5}
Let A be an algebra in $\sC$ and $N$ a module in $\sM(A)$.
Then for each finite family $M_i;\;i\in I$ of modules in
$\sM(A)$, there is a given $\KK$-submodule
$$\Mult_{\sM(A)}(M_i,i\in I;N)$$ of
the module $\Mult_A(M_i,i\in I;N)$ of $A$-multilinear forms
$\prod_{i\in I}M_i\nach N$, which is functorial in each $M_i$ and $N$
and has the following properties:
\begin{enumerate}
\item
Let $i_0$ be an element of $I$ and $u:M'_{i_0}\nach M_{i_0}$
a homomorphism in $\sM(A)$. Set $M''_{i_0}:=\Kokern(u)$
and $M'_i:=M''_i:=M_i$ for $i\in I\setminus\{i_0\}$.
Then the sequence
$$\textstyle{0\rightarrow\Mult_{\sM(A)}(M''_i,i\in I;N)\rightarrow
\Mult_{\sM(A)}(M_i,i\in I;N)\rightarrow\Mult_{\sM(A)}(M'_i,i\in I;N)}$$
induced by $u$ is exact.
\item\label{6252} 
For modules $M,N\in\sM(A)$ there is a canonical isomorphism\\
$\Mult_{\sM(A)}(M;N)\nach\Hom_{\sM(A)}(M;N)$.
\item\label{6253}
For $M$ in $\sM(A)$, the multiplication map $\mu_M:A\times M\nach M$
is in $\Mult_{\sM(A)}(A\times M;M)$.
\item
If $\sigma:I\nach J$ is a bijective map, then the restriction
of the isomorphism
$$\Mult_A(M_i,i\in I;N)\nach\Mult_A(M_{\sigma^{-1}(j)},j\in J;N)$$
defined by $\sigma$, defines an isomorphism
$$\Mult_{\sM(A)}(M_i,i\in I;N)\nach
\Mult_{\sM(A)}(M_{\sigma^{-1}(j)},j\in J;N).$$
\item\label{6255}
Each homomorphism $\phi:A\nach B$ in $\sC$ induces a cartesian
diagram
$$\xymatrix{\Mult_{\sM(B)}(M_i, i\in I;N)\ar[d]\ar[r] &
\Mult_{\sM(A)}((M_i)_{[\phi]}, i\in I;N_{[\phi]})\ar[d]\\
\Mult_B(M_i, i\in I;N)\ar[r]&
\Mult_A((M_i)_{[\phi]}, i\in I;N_{[\phi]})
}$$
\item\label{6256}
For each $i\in I$ let $L_j,\;j\in J_i$ a nonempty finite family of
modules in $\sM(A)$. Set $J:=\amalg_{i\in I}J_i$. Then
the canonical map
$$(\prod_{i\in I}\Mult_{\sM(A)}(L_j,j\in J_i;M_i))\times
\Mult_{\sM(A)}(M_i,i\in I;N)\nach\Mult_A(L_j,j\in J;N)$$
factorises through $\Mult_{\sM(A)}(L_j,j\in J;N)$.
\item
The functor $N\mapsto\Mult_{\sM(A)}(M_i,i\in I;N)$ on $\sM(A)$
is represented by a module ${\bigotimes_{i\in I}}^{\sM}_{A}M_i$ in
$\sM(A)$.
\item
If $I$ is a disjoint union $\cup_{j\in J}I_j$ with $I_j\neq\emptyset$ for
all $j$, then the canonical homomorphism
$${\bigotimes_{i\in I}}^{\sM}_AM_i\nach{\bigotimes_{j\in J}}^{\sM}_A(
{\bigotimes_{i\in I_j}}^{\sM}_A M_i)$$
is an isomorphism in $\sM(A)$.
\item
The canonical map $A\otimes_A^{\sM}M\nach M$ is an isomorphism in $\sM(A)$. 
\end{enumerate}
\item
Let $\phi:A\nach B$ be a homomorphism in $\sC$ and $M$ a module in
$\sM(A)$ and $N$ a module in $\sM(B)$. Then $N\ot_A^{\sM}M$ is
via the canonical $A$-bilinear map\footnote{The existence of this map is a 
consequence of (2), (5.7) and (5.8).}
$$B\times N\ot_A^{\sM}M\nach N\ot_A^{\sM}M$$ a module in $\sM(B)$.
The analogue statement holds for $M\ot_A^{\sM}N$.
\item
Let $k\nach A$ and $k\nach B$ be two homomorphisms in $\sC$ and 
$\phi$ resp. $\psi$ the canonical maps of $A$ resp. $B$ 
in $C:=A\ot_k^{\sM}B$. Let $M$ be a module in $\sM(k)$ and
$\rho:C\times M\nach M$ an element of $\Mult_{\sM(k)}(C\times M;M)$
such that
\begin{enumerate}
\item[(a)]
$\rho$ extends the multiplication of $k$ on $M$.
\item[(b)]
$M$ is via $\rho$ a $C$-module.
\item[(c)]
$M_{[\phi]}$ is in $\sM(A)$ and $M_{[\psi]}$ in $\sM(B)$.
\end{enumerate}
Then $M$ is in $\sM(C)$.
\item
For algebras $A$ and $B$ in $\sC_k$, the canonical map
$A\ot_k^{\sM}B\nach A\ot_k^{\sC}B$ is an isomorphism in $\sM(k)$.
\end{enumerate}

We specify admissible pairs by the following axioms:
\begin{axioms}
Let $A$ be an algebra in $\sC$.
\begin{liste}
\item[(S1)]
When $u:M\nach N$ is a homomorphism of finite $A$-modules in $\sM(A)$,
then the cokernel of $u$ in $\sM(A)$ coincides with the cokernel
of $u$ in $A$-$\Mod$ and for $N=A$ the cokernel of $u$ is
an algebra in $\sC_A$ with respect to the canonical projection
$A\nach\Kokern(u)$.
\item[(S1')]
For any homomorphism $u:M\nach N$ of $A$-modules
the cokernel of $u$ in $\sM(A)$ coincides with the cokernel
of $u$ in $A$-$\Mod$ and for $N=A$ the cokernel of $u$ is
an algebra in $\sC_A$ with respect to the canonical projection
$A\nach\Kokern(u)$.
\item[(S2)]
Bijective homomorphisms in $\sM(A)$ are isomorphisms.
\end{liste}
\end{axioms}
\vspace{-0.6cm}
\begin{beisp}\label{exop}
\begin{liste}
\item
Let $\sC^{(0)}$ be the category of all commutative
$\KK$-algebras and $\sM^{(0)}$ the category of 
modules over algebras in $\sC^{(0)}$. Then
$(\sC^{(0)},\sM^{(0)})$ is an admissible pair of categories that
satisfies axioms (S1') and (S2).
\item
In the first example, we can replace $\sC^{(0)}$ by the category
of all noetherian, commutative $\KK$-algebras.
\item
Let $\sC^{(1)}$ be the category of all 
analytic $\CC$-algebras, i.e. the category of all sections
of the structure sheaf of a Stein compact.
Then each algebra in $\sC^{(1)}$ is a DFN-algebra
and each homomorphism of such algebras is continuous.
Let $\sM^{(1)}$ be the category of all DFN-modules over algebras in
$\sC^{(1)}$. For modules $M$ and $N$ in $\sM^{(1)}$, we set
$\Hom_{\sM^{(1)}}(M,N)$ to be the group of all continuous
homomorphisms $M\nach N$ in $\Hom_{\sC^{(1)}}(M,N)$.
We set $\Mult_{\sM^{(1)}}(\;)$ to be the group of all
continuous multilinear forms. Then $(\sC^{(1)},\sM^{(1)})$ is an 
admissible pair of categories that
satisfies axioms (S1) and (S2).
\item
In the last example, we can replace $\sC^{(1)}$ by the category
of local analytic algebras.
\end{liste}
\end{beisp}

\vspace{-1cm}
\subsection{About the tensor product $\ot_k^{\sM}$}

Let $k$ be an algebra in $\sC$. In the following considerations
$A,B,M$ and $N$ are modules in $\sM(k)$.
By axiom~(5.2), there is a natural isomorphism
$\sim:\Mult_{\sM(k)}(A\times B,M)\nach\Hom_{\sM(k)}(A\ot_k^{\sM}
B,M)$. 
This means that each
morphism $f:M\nach N$ in $\sM(k)$ induces a commutative diagram
\begin{equation*}
\xymatrix{
\Mult_{\sM(k)}(A\times B,M)\ar[r]^{\isom}\ar[d]^{f^{\ast}} &
\Hom_{\sM(k)}(A\ot_k^{\sM}B,M)\ar[d]^{f^{\ast}}\\
\Mult_{\sM(k)}(A\times B,N)\ar[r]^{\isom} & \Hom_{\sM(k)}(A\ot_k^{\sM}B,N)
}
\end{equation*}

We denote the inverse map of $\sim$ also by $\sim{ }$.
\begin{bem}
For $h\in\Hom_{\sM(k)}(A\ot_k^{\sM} B,N)$ we have 
$\tilde{h}=h\circ\tilde{\id}_{A\ot B}$.
\end{bem}
\bew
In the diagram above, set $M:=A\ot_k^{\sM}B$ and $f:=h$.
We have $h=h^{\ast}(\id_{A\ot B})$. So $\tilde{h}$ is the image of
$\id_{A\ot B}$ by going through the diagram starting up right, going
down and left. $h\circ\tilde{\id}_{A\ot B}$ is the result, choosing
the other way.
\qed

The following consequence reminds of the situation, where
$A\ot^{\sM}B$ is a topological tensor product in which the usual
tensor product is a dense subset, and $\Hom_{\sM}$ stands for continuous
mappings.
\begin{kor}
Each homomorphism $h\in\Hom_{\sM(k)}(A\ot_k^{\sM} B,N)$ is uniquely determined
by its values on the elements of the form
$a\ot b=\tilde{\id}_{A\ot B}((a,b))$.
\end{kor}

Now suppose that $A$ and $B$ are $k$-algebras in $\sC$.
We will explain that there are two ways to see the elements
$a\ot b$ in $A\ot_k^{\sM}B=A\ot_k^{\sC}B$:
By the universal property of fibered products, there is a
natural homomorphism of $k$-algebras
$\alpha:A\ot_k^{\Alg}B\nach A\ot_k^{\sC}B$.
\begin{bem}
For elements $a\ot b$ of $A\ot_k^{\Alg}B$, we have
$\alpha(a\ot^{\Mod} b)=\tilde{\id}_{A\ot_k^{\sM}B}((a,b))$.
\end{bem}
\bew
We see that $\alpha$ is just the image of $\id_{A\ot_k^{\sM}B}$
by the composition of the mappings
\begin{align*}
\Hom_{\sM}(A\ot_k^{\sM}B,A\ot_k^{\sM}B)\isom
\Mult_{\sM}(A\times B,A\ot_k^{\sM}B)\incl\\
\Mult_{k-\Mod}(A\times B,A\ot_k^{\sM}B)\isom
\Hom_{k-\Mod}(A\ot_k^{\Mod}B,A\ot_k^{\sM}B).
\end{align*}
\qed
\begin{kor}
For elements $a,a'\in A$ and $b,b'\in B$, we have
$(a\ot b)(a'\ot b')=aa'\ot bb'$. 
\end{kor}

Remark that in the antisymmetrical graded
context (see below) we will have
$(a\ot b)(a'\ot b')=(-1)^{ba'}aa'\ot bb'$
for homogeneous $a,b,a',b'$.

\subsection{Free modules and algebras}\label{fmaa}

In this subsection we remind the definitions of free objects in
the categories $\sC$ and $\sM$, where $(\sC,\sM)$ is an admissible pair of
categories. 
\paragraph{Free algebras:}
A \textbf{marking} on $\sC$
is a family $(F_\tau)_{\tau\in\T}$ of subfunctors
$F_\tau:\sC\nach(\text{sets})$ of the identity functor, such that
$F_\tau(A)$ contains $0$ for all $\tau$ and all objects $A$ in $\sC$.
For a given object $k$ of $\sC$ and a family $(\tau_i)_{i\in I}$, we
consider the functor $F_{I,k}:A\mapsto\prod_{i\in I}F_{\tau_i}(A)$
on the category $\sC_k$.
If $F_{I,k}$ is representable, i.e. there is a $k$-algebra $A$
and a canonical bijection
$$b:\Hom^{\sC}_k(A,B)\nach\prod_{i\in I}F_{\tau_i}(B)$$
for each algebra $B$ in $\sC_k$, then  $A$ together with the family
$(e_i)_{i\in I}=b(\id_A)$ is called the \textbf{free algebra} over
$k$ with free algebra generators $e_i,\;i\in I$. We will write
$A=k\lan e_i\ran_{i\in I}$.
The marking $F$ is called \textbf{representable}, if $F_{I,k}$
is representable for each $k$ in $\sC$ and each finite
family $(\tau_i)_{i\in I}$.

\paragraph{Free modules:}
A \textbf{marking} on $\sM$ is a family $(G_u)_{u\in U}$ of subfunctors
$G_u:\sM\nach(\text{sets})$ of the identity functor, such that for
each $u\in U$ the following condition holds:
For each homomorphism $\phi:A\nach B$ in $\sC$ and each module $N$
in $\sM(B)$ we have $G_u(N_{[\phi]})=G_u(N)$.
For a given algebra $A$ in $\sC$ and a family $(u_i)_{i\in I}$,
we consider the functor $G_{I,A}:M\mapsto\prod_{i\in I}G_{u_i}(M)$
on the category $\sM(A)$.
If $G_{I,A}$ is representable, i.e. there is an $A$-module $M$ and a
canonical bijection $$b:\Hom_{\sM(A)}(M,N)\nach\prod_{i\in I}G_{u_i}(N)$$
for each $A$-module $N$, then $M$ together with the family
$(e_i)_{i\in I}=b(\id_M)$ is called the \textbf{free module} over $A$
with free module generators $e_i,\;i\in I$. We will write
$M=\coprod_{i\in I}Ae_i$.
The marking $G$ is called \textbf{representable}, if $G_{I,A}$ is
representable for each $A$ in $\sC$ and each finite family
$(u_i)_{i\in I}$.\\

\textbf{A marking} on $(\sC,\sM)$ is a pair $(F,G)$ of a marking
$F=(F_\tau)_{\tau\in\T}$ on
$\sC$ and a marking $G=(G_u)_{u\in U}$ on $\sM$ together with a map
$\eta:\T\nach U$, such that $F_{\tau}(A)\sub G_{\eta(\tau)}(A)$
for each $A$ in $\sC$ and each $\tau$ in $\T$.

\begin{axioms}
Let $(F,G)$ be a marking on $(\sC,\sM)$.
\begin{enumerate}
\item[(F1)]
$F$ is representable.
\item[(F2)]
Let $k$ be an algebra in $\sC$ and $A=k\lan e_i \ran_{i\in I}$
be a free $k$-algebra in $\sC$. Then the canonical
homomorphism $k[e_i]_{i\in I}\nach k\lan e_i \ran_{i\in I}$ 
in $k$-$\Mod$ is flat and
the functor $M\mapsto A\ot_k^{\sM}M$ is exact on the category
of finite modules in $\sM(k)$.
\item[(F3)]
Let $A$ be like in (F2) and $A'=k\lan e'_i\ran_{i\in I'}$
be another free $k$-algebra in $\sC$ with $I\sub I'$.
Then $A'$ is flat over $A$ via the homomorphism $A\nach A'$
with $e_i\mapsto e'_i$.
\item[(F4)]
$G$ is representable.
\item[(F5)]
For each $u\in U$ and each $A$ in $\sC$,
$G_u$ is a right exact functor on $\sM(A)$.
\item[(F6)]
Let $A$ be an algebra in $\sC$ and $E=\coprod_{i\in I}Ae_i$ be
a free $A$-module with respect to $G$ with finite basis $(e_i)_{i\in I}$
and let $M$ be a module in $\sM(A)$.
Then the canonical homomorphism $M^I\nach M\ot_AE$ in $A$-$\Mod$
is bijective.
\item[(F7)]
Let $k$ be an algebra in $\sC$ and $A=k\lan e_i \ran_{i\in I}$
be a free $k$-algebra in $\sC$ with finite $I$.
Then $\Omega_{A/k}$ is a free $A$-module with basis
$de_i\in G_{\eta(\tau_i)}(\Omega_{A/k})$.
\end{enumerate}
\end{axioms}
Remark that
Axiom (F2) implies that free algebra generators (of degree 0)
are no zero divisors.

\begin{defi}
The marking $(F,G)$ is called \textbf{good}, if axioms (F1), (F4),
(F5), (F6) and (F7) hold.
An admissible pair of categories $(\sC,\sM)$ equipped with a good
marking $(F,G)$ is called a \textbf{good pair of categories}, if
it satisfies axioms (S1) and (S2).
\end{defi}
\begin{beisp}\label{makke}
\begin{liste}
\item
On the admissible pair $(\sC^{(0)},\sM^{(0)})$ 
of example~\ref{Apoc}.\ref{exop}, we work
with the \textbf{trivial marking}, i.e.  $F(A)=A$ for
each algebra $A$ in $\sC$, $G(M)=M$ for each module $M$ in
$\sM(A)$. With this marking, $(\sC^{(0)},\sM^{(0)})$
is a good pair of categories, that satisfies 
additionally axioms (F2) and (F3).
\item
Consider the admissible pair $(\sC^{(1)},\sM^{(1)})$.
For $A$ in $\sC^{(1)}$ and $t\in T:=(0,\infty)$
let $F_t(A)$ be the set of all elements of $A$,
such that the transformation of Gelfand (see \cite{BAlg} for the
definition)
$\chi(A)\nach\CC$ factorises through $\{z\in\CC:\;|z|\leq t\}$.
Further, let $G$ be the canonical marking on $\sM^{(1)}$.
Then the pair $(\sC^{(1)},\sM^{(1)})$, together with the marking
$(F,G)$ is a good pair of categories, that
satisfies axioms (F2) and (F3).
\item
When $\sC$ is the category of local analytic algebras
and $\sM$ the category of DFN-modules over $\sC$,
then for $G$ we use the trivial marking
and for objects $A$, we set $F(A)$ to be
the maximal ideal $\m_A$ of $A$. 
Then $(\sC,\sM)$ is a good pair of categories, that also
satisfies axioms (F2) and (F3).
\end{liste}
\end{beisp}

\vspace{-1cm}
\subsection{About graded objects}\label{ago}

Let $(\sC,\sM)$ be an admissible pair of categories.
As in \cite{BinKos}, we can construct a new admissible pair
$(\grc,\grm)$ as follows:\\
Let $\grc$ be the category of graded anti-commutative\footnote{i.e.
$ab=(-1)^{g(a)g(b)}ba$ for homogeneous $a,b\in A$} rings
$A=\coprod_{i\leq 0}A^i$ with $A^0$ in $\sC$, all $A^i$ in $\sM(A^0)$,
such that the multiplication maps
$A^i\times A^j\nach A^{i+j}$ belong to 
$\Mult_{\sM(A^0)}(A^i\times A^j,A^{i+j})$.
A homomorphism $\phi:A\nach B$ in $\grc$ is a homomorphism
of graded anti-commutative rings, such that $\phi^0$
is a homomorphism in $\sC$ and all $\phi^i:A^i\nach B^i$ are $\phi^0$-
linear homomorphisms in $\sM$.\\

Let $\grm$ be the category over $\grc$ whose objects
over an algebra $A$ in $\grc$ are the graded $A$-modules
$M=\coprod_{i\in\ZZ} M^i$, with $M^i=0$ for almost all $i>0$,
such that each $M^i$ is in $\sM(A^0)$ and the maps
$A^i\times M^j\nach M^{i+j}$ belong to
$\Mult_{\sM(A^0)}(A^i\times M^j, M^{i+j})$. When $B$ is another algebra
in $\grc$ and $N$ a module in $\grm(B)$, then
$\Hom_{\grm}(M,N)$ is the set of all pairs
$(\phi,f)$, where $\phi:A\nach B$ is a homomorphism in $\grm$
and $f:M\nach N$ is a $\phi$-linear homomorphism of degree zero,
such that all $f^i:M^i\nach N^i$ are homomorphisms in $\sM$
over $\phi^0$.\\

For modules $M_1,\ldots,M_n$ and $N$ in $\grm(A)$, let
$\Mult_{\grm(A)}(M_1\times...\times M_n,N)$ be the $\KK$-module
of all maps $f:M_1\times\ldots\times M_n\nach N$ with the following
properties:
\begin{enumerate}
\item
For $k_1,\ldots,k_n$ in $\ZZ$, the restriction
$f|_{M_1^{k_1}\times\ldots\times M_n^{k_n}}$ factorises
over a map in\\ $\Mult_{\sM(A^0)}
(M_1^{k_1}\times\ldots\times M_n^{k_n},N^{k_1+\ldots+k_n})$.
\item\label{twoo}
For elements $a\in A$ and $m_\mu\in M_\mu$, we have\\
$f(m_1,\ldots,m_ra,m_{r+1}\ldots m_n)=f(m_1,\ldots,m_r,am_{r+1},\ldots,m_n)$
for $1\leq r<n$ and\\
$f(m_1,\ldots,m_na)=f(m_1,\ldots,m_n)a$
\end{enumerate}

We made use of the fact that we can make each
$M$ in $\grm(A)$ an anti-symmetrical $A$-bimodule by
setting $m\cdot a:=(-1)^{g(m)g(a)}a\cdot m$ for homogeneous
elements $a\in A$ and $m\in M$.\\

To define \textbf{free algebras in} $\grc$, we have to modify the
definition in subsection~\ref{fmaa} as follows:
There is a map $g:\T\nach\ZZ_{\leq 0}$ and each functor $F_\tau$ is
a subfunctor of the functor $A\mapsto A^{g(\tau)}$.
In this graded context, when $F$ is representable, then each functor
$F_{I,A}$ is representable, when for $n\leq 0$, the set of all
$\tau_i$ with $g(\tau_i)=n$ is finite. In this case the
free algebra $A\lan e_i\ran_{i\in I}$ is called
\textbf{g-finite} free algebra. Of course, the degree
$g(e_i)$ of a free generator $e_i$ is just $g(\tau_i)$.\\

To define free modules in $\grm$, we have to modify the
definition in subsection~\ref{fmaa} as follows:
There is a map $g:U\nach\ZZ$ and each functor $G_u$ is
a subfunctor of the functor $M\mapsto M^{g(u)}$.
In this graded context, when $G$ is representable, then each functor
$G_{I,A}$ is representable, when for each $n$, the set of all
$u_i$ with $g(u_i)=n$ is finite. In this case the
free module $\coprod_{i\in I}A e_i$ is called
\textbf{g-finite} free $A$-module. We have $g(e_i)=g(u_i)$ for $i\in I$.\\
 
To define a marking on $(\grc,\grm)$, we have to add in 
subsection~\ref{fmaa} the condition, that the map $\eta:\T\nach U$
is compatible with $g$.\\

\begin{beisp}
When $G$ is a marking on $\sM$, there is a marking 
$\gr(G)=(\gr(G))_{u'\in U'}$ on
$\grm$, defined in the following way:
Set $U':=U\times\ZZ$. For $A$ in $\grc$, $M$ in $\grm$ and
$u'=(u,n)\in U'$ set $\gr_{u'}(G)(M):=G_u(M^n)$. 
Here we have $g(u')=n$.\\

When $(F,G)$ is a marking on $(\sC,\sM)$, there is a marking
$\gr_G(F)=(\gr_G(F)_{\tau'})_{\tau'\in\T'}$, defined in the
following way:
Let $\T'$ be the disjoint union of $\T\times\{0\}$ and
$U\times\ZZ_{<0}$. For $A$ in $\grc$ and $\tau'=(\tau,n)$ we set
$\gr_G(F)_{\tau'}(A)=F_\tau(A)$ if $n=0$ and 
$\gr_G(F)_{\tau'}(A)=G_\tau(A^n)$ if $n<0$.\\

When $(F,G)$ is a marking on $(\sC,\sM)$, then
$(\gr_G(F),\gr(G))$ is a marking on $(\grc,\grm)$ with
the map $\eta':\T'\nach U'$ given by
$(\tau,0)\mapsto(\eta(\tau),0)$ and
$(u,n)\mapsto(u,n)$ for $n<0$.\\
\end{beisp}

Remark that by \cite{BinKos}, lemma (7.6),  
free algebra generators
of negative degree behave much like
polynomial variables\footnote{For a more precise statement
see proposition~\ref{Apoc}.\ref{polygod}.}
and when $(\sC,\sM)$ is a good pair of 
categories, then $(\grc,\grm)$ is also a good pair of
categories.  

When $(\sC,\sM)$ is an admissible pair of categories that satisfies axiom (S2),
then by \cite{BinKos} proposition (6.9), the admissible pair 
$(\gr(\sC),\gr(\sM))$ also
satisfies axiom (S2). In general this is not true for axiom (S1). But we have:

\begin{bem}\label{318021}
Let $(\sC,\sM)$ satisfy axiom (S1) and let $A$ be an object of
$\gr(\sC)$, such that all $A^i$ are finite $A^0$-modules. 
Then for g-finite modules $M,N$ in $\gr(\sM)(A)$
each homomorphism $f:M\nach N$ in $\gr(\sM)(A)$ is strict, i.e.
the cokernel of $f$ in $\grm$ coincides with the set-theoretical
cokernel.
\end{bem}
\begin{bem}\label{kopp}
Let $(\sC,\sM)$ be an admissible pair of categories with a marking 
$(F,G)$, where $G$ is trivial. Suppose that axiom $(S1)$ holds.
Let $k$ be an algebra in $\sC$ and let $M_1,M_2$ and $N$ be modules
in $\sM(k)$, such that $M_1$ and $M_2$ are finite $k$-modules with
$M_i\sub N$ and $M_1\cap M_2=\{0\}$. Then we have
\begin{enumerate}
\item
The inclusions $M_i\incl N$ are homomorphisms in $\sM(k)$.
\item
$M_1+M_2$ is in $\sM(k)$.
\item
The inclusions $M_i\nach M_1+M_2$ are homomorphisms in $\sM(k)$.
\item
The projections $p_i:M_1+M_2\nach M_i$ are homomorphisms in $\sM(k)$.
\item
$M_1+M_2=M_1\oplus M_2$. 
\end{enumerate}  
In $(\grc,\grm)$ the same statement is true when we suppose that
all $k^i$ are finite $k^0$-modules and $M_1,M_2$ and $N$ are g-finite.
\end{bem}
\bew
(1) There are free finite $k$-modules $L_i$ in $\sM(k)$ and
homomorphisms $\phi_i:L_i\nach N$ such that the inclusions 
$M_i\incl N$ is $\ker(\kokern(\phi_i))$. (2)
We have $M_1+M_2\incl N=\ker(\kokern(\phi_1+\phi_2))$.
(3) $M_i\incl N$ factorises through $\ker(\kokern(\phi_1+\phi_2))$.
(4) The projection $M_1+M_2\nach M_1$ is the kernel
of the inclusion $M_2\nach M_1+M_2$ in $k$-$\Mod$, so as well in
$\sM(k)$. 
(5) Consider homomorphisms $f_i:M_i\nach P$ in $\sM(k)$.
We define a homomorphism $M+N\nach P$ as $f_1\circ p_1+f_2\circ p_2$.
Then the diagram
$$\xymatrix{
 M_1+M_2\ar[dr] & M_1\ar[l]\ar[d]^{f_1}\\
M_2\ar[u]\ar[r]^{f_2} & P
}$$
in $\sM(k)$ commutes. The graded case follows in the same manner.
\qed

\begin{bem}\label{axt}
\begin{liste}
\item
Suppose that $(\sC,\sM)$ is a good pair of categories 
and that $k$ is an algebra in $\grc$, such that each $k^i$
is a finite $k^0$-module.
Let $R=k\lan T\ran$ be a free g-finite algebra over $k$ in $\grc$.
Then there is a decomposition
$$R=k\oplus\sum_{t\in T}Rt$$
in the category $\grm(k)$.
\item
Suppose that additionally the marking $G$ on $\sM$ is trivial
and that axiom (F2) holds.
Then, for each $n\geq 0$, there is a decomposition
$$R=k\oplus\sum_{t_1\in T}t_1k\oplus...\oplus
\sum_{t_1,\ldots,t_n\in T}t_1\cdot\ldots\cdot t_n k\oplus
\sum_{t_1,\ldots,t_{n+1}\in T}t_1\cdot\ldots\cdot t_{n+1} R.$$
\end{liste}
\end{bem}
\bew
(i)
We can form the free g-finite $R$-module $M=\coprod_{t\in T}Re(t)$,
where to each free algebra generator $t\in\gr_G(F)_{\tau'}$
we have associated a free module generator
$e(t)\in\gr(G)_{\eta(\tau')}(M)$. Now we consider the
homomorphism $M\nach R$ in $\grm(R)$ with $e(t)\mapsto t$.
By remark~\ref{Apoc}.\ref{318021}, 
the cokernel map of this homomorphism coincides with
the cokernel map in $R$-$\Mod$, which is just the projection
$p:R\nach R/(T)$ and $R/(T)$
is an algebra in $\grc$. 
Now there is a diagram
$$\xymatrix{
R\ar[r]^{\pi}\ar[dr]^p&k\ar[d]\\
 & R/(T) 
}$$
in $\grc$, where $\pi:R\nach k$ is the homomorphism
given by $t\mapsto 0$ for $t\in T$ and the homomorphism
$k\nach R/(T)$ is the canonical inclusion. The diagram
commutes, since in both directions a $t\in T$ goes to $0$.
So we get $\Kern(\pi)=(T)$. But obviously, we have
$R=k\oplus\Kern(\pi)$.\\

(2)
The case $n=0$ is just (i).
Instead of doing an induction step, we show the case $n=1$:
Iterating the decomposition in the axiom, we get
$$R=k\oplus\sum_{t_1\in T}t_1k+\sum_{t_1,t_2\in T}t_1t_2R.$$
We have to show that the sum on the right is a direct sum.
With remark~\ref{Apoc}.\ref{kopp}, it is enough to show that
$\sum_{t_1\in T}t_1k\cap\sum_{t_1,t_2\in T}t_1t_2R=\{0\}$.
We can restrict ourselves to the case where $T$ consists
of a single element $t$. If the intersection was not trivial,
then $t$ would be a zero divisor. In the case $g(t)=0$, this  
contradicts (F2). In the case $g(t)<0$ odd, the annulator of
$t$ is generated by $t$. In the case $g(t)<0$ even, $t$ is no zero divisor.
If the intersection was not trivial,
we would find an annulator of $t$, 
which is not in $(t)$. Contradiction!
\qed

\subsection{Simplicial objects}\label{sims}
Let $I$ be an index set. Then a set $\sN$ of subsets of $I$
is called \textbf{simplicial scheme} over $I$,
if $\emptyset\not\in \sN$; if
for all $i\in I$, we have $\{i\}\in\sN$ and if every nonempty
subset of an element in $\sN$ is again in $\sN$.\\

For an element $\alpha$ of a simplicial scheme $\sN$ over $I$,
containing n elements, set $|\alpha|:=n-1$.
Then for $n\geq 0$, the set $\sN^{(n)}$ of
all $\alpha\in\sN$ with $|\alpha|\leq n$ is again a simplicial scheme
over $I$.\\

A simplicial scheme $\sN$ can be seen as category, where
$\Hom(\alpha,\beta)$ contains only the inclusion $\alpha\sub\beta$, if
$\alpha\sub\beta$ and is empty in all other cases.

When $\sA$ is any category, by an $\sN$-object in $\sA$, we mean a covariant
functor $\sN\nach \sA$. The $\sN$-objects in $\sA$ again form a category,
denoted by $\sA^{\sN}$.
When $(\sC,\sN)$ is an admissible pair of categories and
$A=(A_\alpha)_{\alpha\in\sN}$ an object of $\sC^{\sN}$,
we denote the category of $\sN$-objects $M=(M_\alpha)_{\alpha\in\sN}$
in $\sM^{\sN}$ with $M_\alpha\in\ob(\sM(A_\alpha))$ by $\sM^{\sN}(A)$.\\

Let $(\sC,\sM)$ be an admissible pair of categories and
$\sN$ a simplicial scheme. Let 
$((F_\tau)_{\tau\in\T},(G_u)_{u\in U})$ be a marking on $(\sC,\sM)$.
Then for each pair $(\alpha,\tau)$ in $\sN\times\T$, there is a
functor $F_{\alpha,\tau}:A\mapsto F_\tau(A_\alpha)$.
For a family $(\alpha_i,\tau_i),i\in I$ of elements of
$\sN\times\T$ and $A$ in $\sC^{\sN}$, there is a set-valued functor
$B\mapsto\prod_{i\in I}F_{\alpha_i,\tau_i}(B)$. We will denote it by
$F_{I,A}$.

\begin{bem}
Suppose that for each $\alpha\in\sN$ the free $A_\alpha$-algebra
$A_\alpha'=A_\alpha\lan e_i^{(\alpha)}\ran_{\alpha_i\sub\alpha}$
in the free generators $e_i^{(\alpha)}\in F_{\tau_i}(A'_\alpha)$
exists. For $\alpha\sub\beta$ let 
$\rho_{\alpha\beta}:A'_\alpha\nach A'_\beta$ be the
homomorphism in $\sC$ over
$A_\alpha$, given by $e_i^{(\alpha)}\mapsto e_i^{(\beta)}$.
Then $A'=(A'_\alpha)_{\alpha\in\sN}$ is an algebra in $\sC^{\sN}$
and together with the family $(e^{(\alpha_i)})_{i\in I}$,
it represents the functor $F_{I,A}$. We call it the
\textbf{free $A$-algebra} in the free generators
$e_i:=e_i^{(\alpha_i)}\in F_{\alpha_i,\tau_i}(A')$ and denote it
by $A\lan e_i\ran_{i\in I}$.
\end{bem}

For each pair $(\alpha,u)\in\sN\times U$, there is a functor
$G_{\alpha,u}:M\mapsto G_u(M_\alpha)$.
For a family $(\alpha_i,u_i),i\in I$ of elements of
$\sN\times U$ and $A$ in $\sC^{\sN}$, there is a set-valued functor
$N\mapsto\prod_{i\in I}G_{\alpha_i,u_i}(M)$. We will denote it by
$G_{I,A}$.

\begin{bem}
Fix a family $(\alpha_i,u_i),i\in I$ of elements of
$\sN\times U$ and an algebra $A$ in $\sC^{\sN}$
Suppose that for each $\alpha\in\sN$ the free $A_\alpha$-module
$M_\alpha=\coprod_{\alpha_i\sub\alpha}A_\alpha e_i^{(\alpha)}$
in the free generators $e_i^{(\alpha)}\in G_{u_i}(M_\alpha)$
exists. For $\alpha\sub\beta$, let 
$\rho_{\alpha\beta}:M_\alpha\nach M_\beta$ be the
homomorphism in $\sM$ over
$A_\alpha$, given by $e_i^{(\alpha)}\mapsto e_i^{(\beta)}$.
Then $M=(M_\alpha)_{\alpha\in\sN}$ is a module in $\sM^{\sN}$
and together with the family $(e^{(\alpha_i)})_{i\in I}$
it represents the functor $G_{I,A}$. We call it the
\textbf{free $A$-module} with free generators
$e_i:=e_i^{(\alpha_i)}\in G_{\alpha_i,u_i}(A')$ and denote it
by $\coprod_{i\in I} A e_i$.
\end{bem}

To distinguish the nonsimplicial from the simplicial context, we
call the first one \textbf{affine}.

\subsection{Resolutions}\label{reso}

For a DG-module $K$ in $\gr(\sM)(R)$ with differential $d$ -for the
instant- denote
the image of $d^{i-1}:K^{i-1}\nach K^{i}$ in $\sM(R)$ by 
$B^i_{\sM}$ and the kernel of $d^i$ in $\sM(R)$ by $Z^i_{\sM}$.
We use the same letters without subscript to denote image and kernel in
$R$-$\Mod$.
In general, the quotient $Z^i_{\sM}/B^i_{\sM}$
(formed in the category $R$-$\Mod$)
is not an object of $\sM(R)$. So we define the \textbf{i-th homology
$H^i_{\sM}(K)$ of $K$ in} $\sM(R)$ to be the cokernel of the map 
$B^i_{\sM}\nach Z^i_{\sM}$.
When $K$ is separated in the sense that the cokernels of
the maps $d:K^i\nach K^{i+1}$, the induced maps $K^i\nach
Z^{i+1}_{\sM}(K)$ and the inclusions $Z^{i}_{\sM}(K)\nach K^i$
coincide with their cokernels, formed in the category $R$-$\Mod$,
then
$H^i_{\sM}(K)$ is as $R$-module isomorphic to $H^i(K)=Z^i/B^i$.
We call $K$ \textbf{acyclic}, if $H^i(K)=0$ for all i. We call $K$
$\sM$-\textbf{acyclic}, if $H^i_{\sM}(K)=0$ for all i.

\begin{bem}\label{2506021}
Suppose that
$(\sC,\sM)$ satisfies axiom (S1) and 
all $K^i$ are  finite $R$-modules.
Then $K$ is acyclic if and only if $K$ is $\sM$-acyclic.
\end{bem}

By a \textbf{DG-resolution} of an object $B$ in $\sM$,
we mean a DG-module $M$ in $\gr(\sM)$, such that 
$H^i_{\sM}(M)=0$ for $i<0$ and $H^0_{\sM}(M)=B$.

\begin{kor}\label{2506022}
If the pair $(\sC,\sM)$ satisfies (S1) and (S2), then for a DG $R$-module 
$K$ in $\gr(\sM)$ which is finite over $R^0$ in each degree,
the following statements are equivalent:
\begin{enumerate}
\item
$K$ is a DG-resolution of an object $M\in\ob(\sM)$. 
\item
$K$ is a resolution of $M$ as differential graded $R$-module.  
\end{enumerate}
\end{kor}

\vspace{-0.8cm}

We remind definition (8.1) in \cite{BinKos}. Here we work in an
admissible
pair $(\sC,\sM)$ with a fixed marking $(F,G)$.
\begin{defi}
Let $A\nach B$ be a homomorphism of DG-objects in $\grcn$. When we talk
of a \textbf{resolvent of $B$ over $A$}, we mean a free 
DG-algebra $R$ over $A$ in $\grcn$ (with respect to the marking $\gr_G(F)$)
together with a morphism $R\nach B$ of DG-Objects in $\grcn$ which is a
surjective quasi-isomorphism on each $\alpha\in\sN$.
\end{defi}

In this paper, we will mostly work in a noetherian context,
i.e. we will mostly assume that the following axiom is satisfied:
 
\begin{axiom}(N)
Each algebra $A$ in $\sC$ is noetherian and each finite module
$M$ in $\sM(A)$ is a quotient of a finite free $A$-module. 
\end{axiom}

If the good pair $(\sC,\sM)$ satisfies the axioms (N) and (F2) and
if $A^i$ is a finite $A^0$-module for all i and
if $B^i$ is a finite $B^0$-module for all i and if $A^0$ is a quotient
of a g-finite $B^0$-module $C$ in $\sC^{\sN}$, such that
each $C_\alpha$ is a finite free $B_\alpha^0$-algebra, then
such resolvents exist by \cite{BinKos}, prop. (8.7) and prop. (8.8).
We can also deduce
easily those results from remark~\ref{Apoc}.\ref{fres}.
For the existence of resolvents in the non-noetherian case, see 
\textit{loc. cit.}

The next proposition is of great importance for this work.
Here we consider a good pair $(\sC,\sM)$ 
of categories and suppose that the marking $G$ on $\sM$ is trivial
and that axiom (N) is satisfied.

\begin{prop}\label{sup26}
Let $A\nach B$ be a homomorphism of DG-objects in $\grc$.
Then for 
two g-finite resolvents $R_1$ and $R_2$ of $B$ over $A$, 
there exist homomorphisms $R_1\nach R_2$ and $R_2\nach R_1$ in $\grc$,
that
are homotopy equivalences over $A$.
\end{prop}
\bew
\textbf{First case:} Suppose that $R^0_1=R^0_2$.\\
Set $A':=A\ot_{A^0}R_1^0$. Then $R_1$ and $R_2$ are resolvents
of $B$ over $A'$. With \cite{BinKos}, prop. 8.2., there are quasi-isomorphisms
$R_1\nach R_2$ and $R_2\nach R_1$ in $\grc$ over $A'$.
But since $R_1^0=R_2^0={A'}^0$, $R_1$ and $R_2$ are free $A'$-modules
in $\grm$. Hence the quasi-isomorphisms are already homotopy-equivalences.
\\
\textbf{Second case:} Suppose that $R^0_2$ is a finite free algebra over
$R_1^0$ in $\sC$.\\
By induction, we can restrict ourselves to the case, where
$R_2^0=R_1^0\lan e\ran$ is just a free algebra in one generator.
Consider the free $R_1^0$-algebra $R:=R_1^0\lan e,f\ran$,
in $\grc$
generated by a free generator $e$ of degree 0 and a free generator $f$
of degree $-1$. We define a differential on $R$ by setting
$f\mapsto e$. By remark~\ref{Apoc}.\ref{axt} (1), we have 
$R_1^0\lan e\ran=R_0\oplus eR_1^0\lan e\ran$. So by axiom (S2), the
differential gives an isomorphism $fR_1^0\lan e\ran\nach eR_1^0\lan e\ran$.
With this in mind, we can easily construct a contracting homotopy
on $R$. Now $R'_1:=R_1\ot_{R_1^0}R$ is homotopic over $R_1$ to $R_1$.
More precisely, the inclusion $R_1\nach R'_1$ and the projection
$R'_1\nach R'$ are homotopy equivalences.
By the first case, there are homotopy-equivalences
$R'_1\nach R_2$ and $R_2\nach R'_1$.\\
\textbf{General case:}
Let $R_3$ be a free g-finite resolution of $B$ over $R_1\ot_AR_2$.
Now $R_3$ is free over $R_1$ and $R_2$ and by the second case,
there are homotopy-equivalences
$R_1\nach R_3\nach R_2$ and $R_2\nach R_3\nach R_1$ in $\grc$.
\qed

In the simplicial case\footnote{Remember that we still assume that the
marking $G$ on $\sM$ is trivial.}, there is one little difference.
A free algebra in $\gr(\sC)^{\sN}$ over an algebra $A$ in $\gr(\sC)^{\sN}$
is not a free module in $\gr(\sM)^{\sN}(A)$, even if all free algebra
generators are of strictly negative degree. The point is, that even
$A$ itself is not free as $A$-module.
But we see that a free algebra over $A$ with free algebra generators of
negative degree is as $A$-module in $\gr(\sM)^{\sN}$ a direct sum
$A\oplus M$ with a free $A$-module $M$.

\newcommand{\cn}{\sC^{\sN}}
\newcommand{\mn}{\sM^{\sN}}
We need two lemmas to prove a simplicial version of
proposition~\ref{Apoc}.
\ref{sup26}. The first one is a simplicial version of the Comparison Theorem
(for the affine case, see \cite{Weib}, theorem 2.2.6).
\begin{lemma}
Let $A$ be a DG algebra in $\grcn$. Let $P=\coprod_{i\in I}Ae_i$ be a free
DG $A$-module in $\grmn$ with a homomorphism $P^0\nach M$
of $A^0$-modules in $\mn(A^0)$.
Let $N$ be an $A^0$-modules in $\mn$ and $Q$ in $\grmn(A)$ a DG-resolution
of $\N$. Let $\phi:M\nach N$ be a $A^0$-homomorphism in $\mn$.
Then there exists a homomorphism $f:P\nach Q$, lifting $\phi$
and it is unique up to a chain homotopy.
\end{lemma}
\bew
The existence of such an $f$ is not hard to prove. But we only make
use of the uniqueness. So we only prove this part here:
Let $f$ and $g$ two DG-homomorphisms, lifting $\phi$.
Inductively we construct families $\{s_\alpha:\;|\alpha|\leq m\}$
of compatible homotopy maps $s_\alpha:P_\alpha\nach Q_\alpha[-1]$
satisfying
$$g^n-f^n=d_Q\circ s_\alpha^n+(-1)^ns_\alpha^{n+1}\circ d_P.$$
Suppose, that the free generator $e_i$ is associated to the pair
$(\alpha_i,z_i)$ with $\alpha_i\in\sN$ and $z_i<0$.\\
For $m=0$ and each $\beta$ in $\sN$ with $|\beta|=0$, we see that
$P_\beta$ is free DG-module in $\grm(A_\alpha)$, and we can construct
$s_\beta^\bullet$ just as in the affine case.\\
Now suppose that $\{s_\alpha:\;|\alpha|\leq m\}$ is already
constructed. Then for each $\beta\in\sN$ with $|\beta|=m+1$
we have $$P_\beta=\coprod_{\alpha_i\sub\beta}A_\beta e_i.$$
For $\alpha\sub\beta$, denote the restriction map $P_\alpha\nach\P_\beta$
by $\rho_{\alpha\beta}$.
For free generators $e_i$ with $\alpha_i\sub\beta$
but $\alpha_i\neq\beta$, set
$s_\beta(e_i):=\rho_{\alpha_i\beta}(s_{\alpha_i}(e_i))$.
Then we get
\begin{align*}
(g_\beta-f_\beta)(e_i)=
\rho_{\alpha_i\beta}((g_{\alpha_i}-f_{\alpha_i})(e_i))=
\rho_{\alpha_i\beta}([s_{\alpha_i},d_{\alpha_i}](e_i))=
[s_\beta,d_\beta](e_i).
\end{align*}
For free algebra generators $e_i$ with $\alpha_i=\beta$
and say $n=z_i=g(e_i)$,
exactly as in the affine case, by induction on $n$,
we can find elements
$m_i$ in $P_\beta^{z_i-1}$ such that
$$(g_\beta-f_\beta)(e_i)=s_\beta(d(e_i))+(-1)^nd(m_i).$$
Then we set $s_\beta(e_i):=m_i$.\\
In this manner, we get a family $(s_\alpha)_{\alpha\in\sN}$ of
compatible chain homotopies.
\qed

\begin{lemma}\label{pr26}
Let $A$ be a DG algebra in $\grcn$
such that each $A^i$ is a finite $A^0$-module.
Let $M=\coprod_{i\in I}Ae_i$ and
$N=\coprod_{j\in J}Ae_j$ two g-finite free DG $A$-modules
in $\grmn$, such that
all generators $e_i$ and $e_j$ are of negative degree.
Suppose that there is a quasi-isomorphism
$$f=\id_A\oplus f':A\oplus M\nach A\oplus N.$$
Then $f$ is already a homotopy equivalence. More precisely,
there is a homomorphism
$$g=\id_A\oplus g':A\oplus N\nach A\oplus N$$ of DG-modules
and a map $s_\ast:M\nach M[-1]$ of graded modules, such that
$s^0_\ast=0$ and $g\circ f-\id=[s,d]$.
\end{lemma}
\bew
Consider the following diagram, where the first line is just the
mapping cone $\cone(f)=N\oplus M[1]$
of $f$ and the vertical maps are the canonical inclusions:
\begin{equation*}
\xymatrix{
...\ar[r]& M^{-1}\oplus N^{-2}\ar[r]& M^{0}\oplus N^{-1}\ar[r] & N^{0}\\
...\ar[r]& N^{-2}\ar[u]^\iota\ar[r]& N^{-1}\ar[u]^\iota
\ar[r]& N^{0}\ar[u]^\iota
}
\end{equation*}
Since $f$ is a quasi-isomorphism, the mapping cone of $f$ is acyclic,
so the first line is a resolution of the module $\{0\}$.
The map $\iota$ of DG-modules is a lifting of the trivial map
$0 \nach 0$. The zero map $N\nach\cone(f)$ is a second candidate for
such a lifting. So we are almost in the situation of the uniqueness
statement in the comparison theorem. The only difference is, that
$N=A\oplus N'$ is not a free module in $\grmn$. But to construct
a chain homotopy $\sigma:N\nach\cone(f)[-1]=N[-1]\oplus M$
for $0\simeq\iota$, we can set $\sigma|_A$ to be the composition
of the inclusions $A\nach A\oplus M'=M$ and $M\nach\cone(f)[-1]$.
On the free generators of $N'$, the map $\sigma$ can be defined
exactly as in the proof of the comparison theorem.
So we can work with a
family of maps
$$\sigma^n=(g^n,t^n):N^n\nach M^{n}\oplus N^{n-1}$$
for $n\leq 0$,
satisfying the condition
$$\iota^n=\delta^{n-1}\sigma^n+(-1)^n\sigma^{n+1}d^n.$$
Here $d$ denotes the differential of $N$ and $\delta$ the differential of
$\cone(f)$.
The evaluation of this conditions shows, that
$g$ is a chain map $N\nach M$ and that $t$ is a chain
homotopy for $\id_N\simeq f\circ g$.\\
In an analogue manner, we get a chain map $h:M\nach N$ with
$\id_M\simeq g\circ h$. We see easily that then we have
$h\simeq f$, so we get $\id_M\simeq g\circ f$.
\qed

Of course, we can also show that two free module resolutions
of a module in $\grmn$ are homotopy equivalent over the
base ring. Now we can state the announced simplicial version of
proposition~\ref{Apoc}.\ref{sup26}:
Suppose that the same assumptions hold.

\begin{prop}\label{rehoeq}
Let $A\nach B$ be a homomorphism of DG-objects in $\grcn$.
Then, for 
two g-finite resolvents $R_1$ and $R_2$ of $B$ over $A$, 
there exist homomorphisms $R_1\nach R_2$ and $R_2\nach R_1$ in $\grcn$,
that
are homotopy equivalences over $A$.\\ 
\end{prop}
\bew
We imitate the proof of prop~\ref{Apoc}.\ref{sup26}.
For the first step we have to use lemma~\ref{Apoc}.\ref{pr26}.
The second and third step are easy to generalize.
\qed

\subsection{Double graded objects}\label{dgo}

Let $(\sC,\sM)$ be an admissible pair of categories.
We define the pair $(\gr^2(\sC),\gr^2(\sM))$ as follows:
The objects of $\gr^2(\sC)$ are the double graded rings
$A=\coprod_{i,j\leq 0}A^{i,j}$ with $A^{0,0}$ in $\sC$ and
all $A^{i,j}$ in $\sM(A^{0,0})$ such that
\begin{enumerate}
\item
For $a\in A^{i,j}$ and $b\in A^{k,l}$ we have
$ab=(-1)^{(i+j)(l+k)}ba$.
\item
The multiplication maps $A^{i,j}\times A^{k,l}\nach A^{i+j,k+k}$
belong to \\
$\Mult_{\sM(A^{0,0})}(A^{i,j}\times A^{k,l},A^{i+j,k+l})$.
\end{enumerate}
Following the ideas of subsection~\ref{ago}, we
can define
$\Hom_{\gr^2(\sC)}(A,B)$ for objects A, B in $\gr^2(\sC)$,
the category $\gr^2(\sM)$, $\Hom_{\gr^2(\sM)}(M,N)$ for
objects $M,N$ of $\gr^2(\sM)$ and $\Mult_{\gr^2(\sM)(A)}(M_1,...,M_n,N)$
for modules $M_1,...,M_n,N$ in $\gr^2(\sM)(A)$.
We don't make this definitions explicit here.

\begin{bem}
Let $A$ be an object of $\gr^2(\sC)$ and $M,N$ objects of
$\gr^2(\sM)$. For $(p,q)$ in $\ZZ\times\ZZ$ set
$T^{p,q}:=\coprod_{i+j=p,k+l=q}M^{i,k}\ot_{A^{0,0}}^{\sM}N^{j,l}$.
Then $T:=\coprod_{p,q}T^{p,q}$ is a tensor product of
$A$ and $B$ in $\gr^2(\sM)(A^{0,0})$. $T$ can be seen in two different
ways as object of $\gr^2(\sM)(A)$. Consider the homomorphism
$u:A\ot_{A^{0,0}}T\nach T$ in $\gr^2(\sM)(A^{0,0})$, sending
$a\ot m\ot n$ to $ma\ot n-m\ot an$. $u$ can be seen in two manners
as homomorphism in $\gr^2(\sM)(A)$. Both of them induce the same
$A$-module structure on $\bar{T}:=\Kokern(u)$.
We see that $\bar{T}$ is a tensor product of $M$ and $N$ in
$\gr^2(\sM)(A)$. 
\end{bem}

\vspace{-0.8cm}

\begin{bem}
The pair $(\gr^2(\sC),\gr(\sM)^2)$ is an admissible pair of
categories.
\end{bem}
\bew Analogue to the proof of (6.9) in \cite{BinKos}.
\qed

\textbf{Convention:}
When we consider an object $K$ of $\gr(\sM)$ as object
of $\gr^2(\sM)$, we set $K^{i,0}=K^{i}$ and $K^{i,j}=0$ for 
$j\neq 0$.\\

In the same manner as above, we can define a marking 
$(\gr^2_G(F),\gr^2(G))$ on the pair $(\gr^2(\sC),\gr^2(\sM))$:\\
Define the index set $\T''$ as 
$\T\times\{0,0\}\cup U\times(\ZZ^{\leq 0}\times\ZZ^{\leq 0})\setminus(0,0)$.
For $\tau''=(\tau,0,0)\in \T''$ and $A\in\gr^2(\sC)$ set
$\gr^2_G(F)_{\tau''}(A):=F_\tau(A^{0,0})$ and for
$\tau''=(u,p,q)$ in $\T''$ with
$(p,q)\neq(0,0)$ set
$\gr^2_G(F)_{\tau''}(A):=G_u(A^{p,q})$.
Define the index set $U''$ as
$U\times\ZZ\times\ZZ$. For $u''=(u,p,q)\in U''$ and $M\in\gr^2(\sM)$
set $\gr(G)_{u''}(M):=G_u(M^{p,q})$.\\

In analogy to \cite{BinKos}, lemma (7.6), we get:
\begin{prop}\label{polygod}
\begin{liste}
\item[(1)]
Let $A$ be an algebra in $\gr^2(\sC)$ and $A'=A\lan e_i\ran_{i\in I}$
a free algebra over $A$, with respect to the marking $\gr^2_G(F)$.
Suppose that the bidegree of each $e_i$ is different to 0. Then
the canonical homomorphism $A[e_i]_{i\in I}\nach A'$ in
$A$-$\Alg$ is bijective.
\item[(2)]
If $(F,G)$ is good, then $(\gr^2_G(F),\gr^2(G))$ is also good.
\end{liste}
\end{prop}
\vspace{-1cm}

\begin{defi}\label{DGbla}
A \textbf{DG-algebra} in $\gr^2(\sC)$ is an algebra $A$ in $\gr^2(\sC)$
equipped with a (vertical) 
$A^{0,0}$-derivation\footnote{This means that
for homogeneous $a,b\in A$ we have $v(ab)=v(a)+(-1)^{a}av(b)$.
In the exponent, by $a$, we mean the total degree of $a$.}
$v:A\nach A$ of bidegree $(0,1)$ with $s^2=0$.\\ 
A \textbf{DDG-Algebra} in $\gr^2(\sC)$ is a DG-algebra $A$ in 
$\gr^2(\sC)$ equipped with a (horizontal) derivation $h$ of
bidegree $(1,0)$, that anti-commutes with $v$, such that $h^2=0$. 
A \textbf{homomorphism} between (D)DG-algebras
is a morphism in $\gr^2(\sC)$ that commutes with the vertical
(and horizontal) differentials.\\

Now let $(A,s)$ be a DG-algebra in $\gr^2(\sC)$.\\
A \textbf{DG-module} in $\gr^2(\sM)(A)$ is a module $M$ in 
$\gr^2(\sM)(A)$
equipped with a (vertical) differential $v$ of bidegree $(0,1)$
such that for $a\in A$ and $m\in M^{k,l}$ the formula
$v(ma)=v(m)a+(-1)^{k+l} mv(a)$ holds.\\
At least let $(A,v,h)$ be a DDG-algebra in $\gr^2(\sC)$.
A \textbf{DDG-module} in $\gr^2(\sM)(A)$ is a DG-module 
$(M,v)\in\gr^2(\sM)(A)$ 
equipped with a horizontal differential $h$, that commutes
with $v$ and such that for $a\in A$ and $m\in M^{k,l}$ the
formula $h(ma)=e(m)a+(-1)^{k+l} mt(a)$ holds.\\
A \textbf{homomorphism} between (D)DG-modules
is a morphism in $\gr^2(\sM)$ that anti-commutes with the vertical
(and horizontal) differentials.\\
\end{defi}

\vspace{-1cm}

\begin{bem}\label{DDGstr}
Let $K=(K,h,v)$ be a DG-algebra in $\gr^2(\sC)$. We consider
a free algebra $K\lan E\ran$ over $K$ in $\gr^2(\sC)$ with a set 
$E=\{e_i:\;i\in I\}$ of free algebra generators 
with $e_i\in \gr^2_G(F)_{\tau''_i}(K\lan E\ran)$ for a certain
$\tau''_i\in \T''$.
For each $i$, if $g(x_i)\neq(0,0)$
choose an element $h_i\in G_{u_i}(K\lan E\ran^{g(x_i)+(1,0)})$
and an element $v_i\in G_{u_i}(K\lan E\ran^{g(x_i)+(0,1)})$, where
$u_i$ is the first component of $\tau_i=(u_i,g(x_i))$.
Then, setting $h(e_i):=h_i$ and $v(e_i):=v_i$, we get an
extension of the horizontal and the vertical derivation
$h$ and $v$ of $K$. This extensions make $K\lan E\ran$ a
DDG-algebra, if and only if for each $i$, we have
\begin{enumerate}
\item
$h(v_i)+v(h_i)=0$ and
\item
$h(h_i)=v(V_i)=0$.   
\end{enumerate}
\end{bem}

\vspace{-0.8cm}
\bew
Inductively,
we can reduce the proof to the case where $E$ consists
of a single element $e$ of bidegree $(p,q)$. In this case, 
it is an easy calculation.\qed
\begin{defi}
By a \textbf{DG-resolution} of an algebra $B$ in $\gr(\sC)$,
we mean a DG-algebra $A$ in $\gr^2(\sC)$, such that for all i the i-th 
row is a DG-module resolution of $B^i$.
By \textbf{DDG-resolution} of a DG-algebra $B$ in $\gr(\sM)$,
we mean a DDG-algebra $A$ in $\gr^2(\sM)$ that is a 
DG-resolution of $B$ such that the map $A^{*,0}\nach B$
is a homomorphism of DG-Algebras in $\grc$.
\end{defi} 

Now when $A\nach B$ is a homomorphism of DG-Algebras in $\grc$,
to get a resolvent $R$ of $B$ over $A$,
it is enough to construct a DDG-resolution $K$
of $B$, that is free over $A$ as object of $\gr^2(\sC)$.
Then we can choose $R$ to be the total complex $\tot(K)$.
So the question of the existence of DDG-resolutions is quite natural.
An answer to this question is given by the following remark:

\begin{bem}\label{fres}
Suppose that for the pair $(\sC,\sM)$ the axioms (N) and (F2) hold.
Let $K=(K,h,v)$ be a DDG-algebra in $\gr^2(\sC)$ and $u:K^{*,0}\nach A$
a homomorphism of DG-algebras in $\grc$.
\begin{enumerate}
\item
When $A^0$ is a quotient of a free $K^{0,0}$-algebra,
then there exists a free DDG-algebra $L=K\lan F\ran$ over $K$, where
$F$ is a g-finite set of generators of bidegree $(k,0);\;k\leq 0$,
and a surjective homomorphism $L^{*,0}\nach A$ over $K^{*,0}$.
\item
For a fixed $p<0$, suppose that we have $u^{p+1}=\kokern(v^{p+1,-1})$.
Then there is a free DDG-Algebra $L=\lan F\cup G\ran$ over $K$
with finite sets $F$ resp. $G$ of generators of bidegree
$(p,-1)$ resp. $(p+1,-1)$, such that we still have 
$u^{p+1}=\kokern(v^{p+1,-1})$ and additionally
$u^{p}=\kokern(v^{p,-1})$ holds.
\item
Fix $p\leq 0$ and $q\leq -1$. Suppose that we have $H^{q+1}(K^{p+1,*})=0$.
Then there is a free DDG-Algebra $L=K\lan F\cup G\ran$ over $K$
with finite sets $F$ resp. $G$ of generators of bidegree
$(p,q)$ resp. $(p+1,q)$, such that we still have
$H^{q+1}(K^{p+1,*})=0$ and additionally $H^{q+1}(K^{p,*})=0$ holds.
\end{enumerate}
\end{bem}

\vspace{-0.7cm}
\bew
(i) is trivial. The proofs of (ii) and (iii) are very similar, so we only
do the proof of (iii).
We choose $G$ such that there is an epimorphism 
$\pi:\amalg_{g\in G}K^{0,0}g\nach\Kern(v^{p+1,q+1})\cap\Kern(h^{p+1,q+1})$.
Set $v(g)=\pi(g)$ and $h(g)=0$.\\
We choose $F$ such that there is an epimorphism 
$\pi':\amalg_{f\in F}K^{0,0}f\nach\Kern(v^{p,q+1})$.
Set $v(f)=\pi'(f)$ and choose $h(f)$ in $\amalg K^{0,0}g$ such that
we get $v(h(f))=-h(v(f))$.
\qed

\begin{defi}\label{dach}
For a g-finite free DG-module $M=\coprod_{i\in I}Af_i$ in $\grm$
with differential $d$
(this construction can be done more generally in $\grmn$),
we can define the
\textbf{exterior algebra} $\wedge_AM$, to be the free DDG-algebra
$A\lan E\ran=A\lan e_i\ran_{i\in I}$ in $\gr^2(\sC)$, where to each $f_i$ 
we associate a free
algebra generator $e_i$ of bidegree $(g(i),-1)$. We wrote $E$ for
the set of all $e_i$.
The vertical differential of $\wedge_AM$ is set to be trivial,
and the horizontal differential $h$ is defined as follows:
Suppose that $d(f_i)=\sum_ja_{ij}f_j$ for a finite family $a_{ij}$ in
$A$. Then set $h(e_i):=\sum_ja_{ij}e_j$.\\
The total complex of $\wedge_AM$ has the structure of a DG-algebra
in $\grc$ and corresponds to the classical definition of an
exterior algebra.\\
In this situation, let $\wedge^j_AM$ be the DG-module in
$\grm$ with $(\wedge^j_AM)^n=A\lan E\ran^{(n,j)}$ for all $j\geq 0$.
\end{defi}

In particular we have $\wedge^0_AM=A$ and $\wedge^1_AM\isom M$
and 
\begin{equation}\label{totto}
\tot(\wedge_AM)=\coprod_{j\geq 0}\wedge^j M[j].
\end{equation}

\subsection{The (cyclic) bar complex}\label{barr}

Let $(\sC,\sM)$ be an admissible pair of categories.
Consider a homomorphism $k\nach A$ of DG-objects in
$\grc$.\\
By the universal property of fibered products, there are given
two maps $A\nach A\ot_k^{\sC}A$; we denote
them in the sequel by $\iota_1$ and $\iota_2$. We denote the 
multiplication map $A\ot_k^{\sC}A\nach A$ by $\mu$. It is
just the uniquely defined homomorphism such that the
diagram

\begin{equation*}
\xymatrix{
k\ar[r]\ar[d] & A\ar[d]^{\iota_1}\ar@/^/[ddr]^{=}\\
A\ar[r]^{\iota_2}\ar@/_/[drr]_{=} & A\ot_{k}^{\grc}A \ar[dr]^{\mu}\\
& & A
}
\end{equation*}

commutes. When $k\nach A$ is a homomorphism of
DG-algebras and $d$ is the differential of $A$, then
$R:=A\ot^{\grc}_kA$ also is a DG-algebra, whose differential is given by
$s=d\ot 1+1\ot d$ and the homomorphisms $\iota_1, \iota_2$ and $\mu$
are morphisms of complexes.\\
  
Let $M$ be a DG $A$-bimodule in $\gr(\sM)$, which is a symmetrical
$k$-bimodule.
Then we can consider $M$ as DG object of $\grm(R)$, where the
scalar multiplication $R\times M\nach M$ satisfies
$(a\ot a',m)\mapsto(-1)^{a'm}ama'$ for homogeneous elements\footnote{
In the exponents we write sometimes just $a$ instead of $g(a)$ for
homogeneous elements. $ab$ then means $g(a)\cdot g(b)$ and not
$g(ab)$, which is just $g(a)+g(b)$.}
$a,a'\in A$ and $m\in M$.    
To see this we have to apply axioms~(5.3),(5.5) and
(5.6). The same axioms must be used to define the
mappings in the next paragraph.

We now define the \textbf{cyclic barcomplex} 
$C_\bullet^{\cycl}(A,M)=(C_\bullet^{\cycl}(A,M),b)$ as well as
the \textbf{(acyclic) barcomplex}
$C_\bullet^{\ba}(A,M)=(C_\bullet^{\ba}(A,M),b')$
of $A$ with values in $M$ as complex of DG-modules in $\grm(A)$.
For $n=0,1,\ldots$ set $C_n^{\cycl}(A,M):=M\ot A^{\ot n}$ and
$C_n^{\ba}(A,M):=M\ot A^{\ot n}\ot A$.\footnote{Here all tensor 
products are formed in the category $\grm(k)$.} 
We can define homomorphisms
\begin{equation*}
d_i:M\ot\underbrace{A\ot...\ot A}_{n}\nach
M\ot\underbrace{A\ot...\ot A}_{n-1},
\end{equation*}
sending
elements $a_0\ot...\ot a_{n}$ to
$a_0\ot...\ot a_i\cdot a_{i+1}\ot...\ot a_{n}$
$(i=0,...,n-1)$ and a homomorphism $d_{n}$, sending
homogeneous elements $a_0\ot...\ot a_{n}$ to
$(-1)^{a_n(a_1+...+a_{n-1})}a_0\cdot a_n\ot a_1\ot...\ot a_{n-1}$.
Those homomorphisms are homomorphisms in $\grm(A)$, when we regard
the tensor products $M\ot A\ot...$ as $A$-modules by left-multiplication
on the first factor. Remark that when $M$ is only a $A$-right module,
we consider it as an antisymmetrical $A$-bimodule by setting
$m\cdot a:=(-1)^{ma}a\cdot m$.\\
Set $b'_{n-1}:=d_0-...+(-1)^{n-1}d_{n-1}$
and $b_{n}:=b+(-1)^{n}d_{n}$.\\

Exactly as in the algebraic case
(see \cite{BinKos}, chapter III, (2.1)), $b$ and $b'$ define
differentials, i.e.
$b^2=0$ and ${b'}^2=0$. They are homomorphisms over
$A$. So $C_\bullet^{\cycl}(A,M)$ and $C_\bullet^{\ba}(A,M)$  
are complexes in $\grm(A)$. 
Remark that $C^{\ba}_{\cdot}(A,M)$ is even a complex
in $\grm(R)$, when we define the $R$-module structure on
$M\ot A^{\ot n}\ot A$ by
$$(a\ot a')\cdot (m\ot\alpha\ot a_{n+1})=
(-1)^{a(a'+m+\alpha)}a'm\ot\alpha\ot a\cdot a_{n+1}$$
for homogeneous elements $a,a',\alpha$ and $m$.
In the sequel we will write $C_\bullet^{\cycl}(A)$ for $C_\bullet^{\cycl}(A,A)$
and  $C_\bullet^{\ba}(A)$ for $C_\bullet^{\ba}(A,A)$

\begin{bem}
In $\grm(k)$ there exist homomorphisms
\begin{equation*}
h_n:\underbrace{A\ot_k^{\sC}...\ot_k^{\sC} A}_{n}\nach 
\underbrace{A\ot_k^{\sC}...\ot_k^{\sC} A}_{n+1}
\end{equation*}
sending elements $a_1\ot...\ot a_n$ to
 $1\ot a_1\ot...\ot a_n$. They define a contracting homotopy
for the bar complex.
\end{bem}

As consequence we see that the 
acyclic bar complex $C^{\ba}(A)$ complex is acyclic.

\begin{bem}
The double complex $C^{\ba}(A)$,
equipped with the $*$-product\footnote{For the definition of the
$*$-product see \cite{BinKos}; we won't use it here.},
is a DDG-algebra in $\gr(\sC)$ over $R:=A\ot_kA$.
So its total complex can serve as a DG-resolution of $A$ over $R$.
\end{bem}
\bew
\cite{BinKos}, chapter III, theorem (2.2)
\qed

\textbf{Attention:} In the analytic case, $\tot(C^{\ba}(A))$
is \textbf{not} a free object in $\grc$.\\
Now we state two well-known relations between
the cyclic and acyclic bar complexes.
For this we consider $R$ as $A$-bimodule via
$a(a_1\ot a_2)=aa_1\ot a_2$ and $(a_1\ot a_2)a=a_1\ot a_2a$.

\begin{prop}
\begin{liste}
\item
There is an isomorphism $C_\bullet^{\cycl}(A,R)\nach C_\bullet^{\ba}(A)$
of complexes in $\grm(A)$, which is in the n-th component given by
\begin{align*}
C_n^{\cycl}(A,R)&\nach C_n^{\ba}(A)\\
(a\ot a')\ot\alpha&\mapsto (-1)^{a(a'+\alpha)}a'\ot\alpha\ot a
\end{align*}
with $\alpha\in A^{\ot n}$.
\item
There is an isomorphism $C_{\cdot}(A,M)\nach M\ot_R C^{\ba}_{\cdot}(A)$
of complexes in $\grm(A)$, where the differential of the second complex is
given by $1\ot b'$. 
In the n-th component it is given by
\begin{align*}
C_\cdot^{\cycl}(A,M)&\nach M\ot_R C^{\ba}_{\cdot}(A)\\
m\ot\alpha&\mapsto m\ot 1\ot\alpha\ot 1.
\end{align*} 
\end{liste}
\end{prop}
In the algebraic Hochschild theory the cyclic bar complex is
often called Hochschild chain complex and the Hochschild cochain
complex is defined in the algebraic literature as 
the complex $C^{\bullet}(A,M)=(C^{\bullet}(A,M),\beta)$ where $C^0(A,M)=M$
and $C^n(A,M)=\Hom_k(A^{\ot n},M)$ for $n=1,2,\ldots$.
The differential $\beta$ is given by:
\begin{align*}
\beta(f)(a_1,\ldots,a_{n+1})=
a_1f(a_2,\ldots,a_{n+1})-f(a_1\cdot a_2,\ldots,a_{n+1})+\ldots\\
+(-1)^{n}f(a_1,\ldots,a_{n}a_{n+1})+(-1)^{n+1}f(a_1,\ldots,a_{n})a_{n+1}.
\end{align*}
\pagebreak
\begin{prop}
\begin{liste}
\item
When $M$ is an anti-symmetrical $A$-bimodule, then
there is an isomorphism of complexes
$$\Hom_k(A^{\ot n},M)\nach \Hom_A(C_n^{\cycl}(A),M),$$
where the differential on the left complex is $\beta$
and the differential on the right complex is the one induced by
the differential $b$ on $C_{\bullet}^{\cycl}(A)$.
\item
There is an isomorphism of complexes
$$\Hom_R(C^{\ba}_\bullet(A),M)\nach\Hom_k(A^{\ot n},M),$$
sending an $f:C^{\ba}_n\nach M$ to the mapping
$a_1\ot\ldots\ot a_n\mapsto f(1\ot a_1\ot\ldots\ot a_n\ot 1)$.
\end{liste}
\end{prop} 

\vspace{-1cm}

\subsection{Regular sequences}\label{hans}

In this section we want to define
a regular sequence for the graded, anticommutative context.
In our definition, the question if a sequence is regular won't depend
on the order of its elements. 
For this we suppose that the ground ring $\KK$ contains
the rational numbers.\\

Here we work with an admissible pair of categories $(\sC,\sM)$,
equipped with a marking\\
$((F_t)_{t\in T},(G_u)_{u\in U})$, that induces the
marking $((\gr_G(F))_{\tau'\in T'},
(\gr(G))_{u'\in U'})$ on $(\grc,\grm)$ and the marking
$((\gr^2_G(F))_{\tau''\in\T''},(\gr^2)_{u''\in U''})$ on
$(\gr^2(\sC),\gr^2(\sM))$.

\begin{defi}
Let $R$ be an algebra in $\grc$. We call a g-finite
subset $X$ of $R$ 
a \textbf{handy sequence} if for each $x$, there is an $u(x)\in U$
such that $x\in\gr(G)_{(u(x),g(x))}(R)=G_{u(x)}(R^{g(x)})$.
When $R=(R,s)$  is a DG-algebra, then a handy sequence $X\sub$
is called \textbf{handy s-sequence} if we have $s(X)\sub(X)$.
For a handy sequence $X\sub R$, let $E$ be a set of free
algebra generators, containing for each $x\in X$ a generator
$e(x)\in\gr^2_G(F)_{(u(x),g(x),-1)}(R\lan E\ran)$
of bidegree $(g(x),-1)$.
Then we call the free DG-algebra
\footnote{in the sense of definition~\ref{Apoc}.\ref{DGbla}} 
$K(X):=R\lan E\ran$ in $\gr^2(\sC)$ over $R$, where the differential
(of bidegree $(0,-1)$) is given by $e(x)\mapsto x$,
the \textbf{Koszul complex} of $X$ over $R$.  
\end{defi}

For practical reasons, when we work with a handy
sequence $X=\{x_i:\;i\in J\}$,
we define an ordering on the index set $J$,
subject to the condition $g(x_i)\leq g(x_j)$ for $i\leq j$.
Remark that for a handy sequence $X\sub R$ and each subset
$Y\sub X$, the quotient\footnote{By ``quotient'', we mean the cokernel
  in $\grm$ of the embedding $(X)\incl R$.} $R/(Y)$ exists in $\grc$.
And when $R$ is a DG-Algebra $(R,s)$ and $X$ is $s$-handy,
then the quotient $R/(X)$ is also a DG-Algebra.

\begin{defi}\label{regS}\textbf{(and Theorem)} 
Suppose that $\QQ\sub\KK$.\\
Let  $X\sub R$ be a handy sequence and let $I$ be the ideal $(X)\sub R$.
Suppose that for each subset $Y\sub X$, we have
$\cap_{n\geq 1} I^nR/(Y)=0$. Then $X$ is called a \textbf{regular sequence}
if one of the following equivalent conditions holds:
\begin{enumerate}
\item[(i)]
Let $T$ be a set of free algebra generators that contains for
each $x\in X$,
 an element $t(x)$ with $g(t(x))=g(x)$. Then the map
$R/I[T]\nach\gr_I(R)=R/I\oplus I/I^2\oplus...$ in
$\gr(\QQ-\Alg)_{R/I}$, sending
$t(x)$ to the class of $x$ in $I/I^2$ is an isomorphism
of (differential) graded $R/I$-algebras. 
\item[(ii)]
For each $x\in X$ and
for each ideal $J\sub R$ that is generated by a subset $Y\sub X$
with $x\not\in Y$ we have:
If $g(x)$ is even, then $x$ is no zero divisor in $R/J$.
If $g(x)$ is odd, then the annulator of $x$ in $R/J$ is just
the ideal, generated by the class of $x$.
\item[(iii)]
The Koszul complex $K(X)$ is a DG-resolution of $R/(X)$ over $R$.
\item[(iv)]
$H^{-1}(K(X))=0$.
\end{enumerate}
\end{defi}
\newcommand{\inn}{\operatorname{in}}
The implication (iii)$\impl$(iv) is trivial.\\
\textbf{Proof of (i)$\impl$(ii)}\\
For an element $r\in R$, let $n(r)$ be the greatest $n$ such
that $r$ is contained in $I^n$ and let $\inn(r)$
be the element represented by $r$ in $I^{n(r)}/I^{n(r)+1}\sub\gr_I(R)$.
Then for elements $r,r'\in R$ we have
\begin{equation}\label{inn}
\inn(rr')=rr'+I^{n(r)+n(r')+1}.
\end{equation}
\textbf{Claim:} A subset $X\sub R$ satisfies condition (ii),
if the subset $\{\inn(x):\;x\in X\}\sub\gr_I(R)$ satisfies
condition (ii).\\

Proof of the claim:
First step:
For $x\in X$, if $g(x)$ is even and $\inn(x)$ is no zero divisor,
then $x$ is no zero divisor. If $g(x)$ is odd and the annulator
of $\inn(x)$ in $\gr_I(R)$ is the ideal, generated by $\inn(x)$,
then the annulator of $x$ in $R$ is the ideal generated by $x$.\\

The even case follows immediately by (\ref{inn}). In the odd case,
let $r$ be in the annulator of $x$, i.e $rx=0$.
Then by (\ref{inn}), we get $\inn(x)\cdot\inn(r)=0$. So by
the assumption, there is an $y_1\in R$, such that
$\inn(r)=\inn(x)\cdot\inn(y_1)$. This implies that
$r_1:=r-xy_1$ is in $I^{n(r)+1}$ and $n(r_1)\geq n(r)+1$.
Since $x^2=0$, we have $r_1x=rx=0$, and in the same way, we
find a $y_2\in R$ with $r_2:=r_1-xy_2\in I^{n(r-1)+1}$.
Inductively, for each $m>n(r)$, we find $y_1,...,y_k$, such
that $r_k:=r-x(y_1+...+y_k)\in I^m$. So $r$ is an
accumulation point of the ideal $(x)$ in the $I$-adic topology of
$R$. But the closure of $(x)$ in this topology is 
just $\cap_{k\geq 0}((x)+I^k)$, which, by the condition  
$\cap_{n\geq 1} I^nR/(Y)=0$, is equal to $(x)$.
So $r$ is an element of $(x)$. The second inclusion is obvious.\\

Second step:
For $x\in X$, 
if weather $g(x)$ is even and $\inn(x)$ is no zero divisor,
or $g(x)$ is odd and the the annulator of $\inn(x)$ in
$\gr_I(R)$ is $(\inn(x))$,
then $(x)\cap I^{n(x)+n}=xI^n$ for each $n\geq 0$.\\

One inclusion and the even case are easy to see.
Suppose that $g(x)$ is odd and that $rx$ is in $I^{n(x)+n}$.
We have to find $r'\in I^n$ such that $xr=xr'$.
If $r\in I^n$, we are done. Else, we have $n(r)<n$
and $\inn(r)\cdot\inn(x)=rx+I^{n(r)+n(x)+1}=0$.
So there is a $y\in R$ such that
$\inn(r)=\in(x)\cdot\inn(y)$. This means that
$r_1:=r-xy$ is in $I^{n(r)+1}$ and we have
$r_1x=rx$. Inductively, we find an $r':=r_{n-n(r)}$,
such that $r'\in I^n$ and $r'x=rx$.\\

As consequence, when we set $\bar{R}:=R/(x)$ and $\bar{I}:=I/(x)$,
we get an isomorphism
$$\gr_I(R)/(\inn(x))\isom\gr_{\bar{I}}(\bar{R}).$$  
We deduce inductively, that for $\bar{R}:=R/(x_1,...,x_{s})$ 
and $\bar{I}:=I/(x_1,...,x_{s})$,
we get an isomorphism
$$\gr_I(R)/(\inn(x_1),...,\inn(x_s))\isom\gr_{\bar{I}}(\bar{R}).$$  

Last step:
When $g(x)$ is even, we have to show, that $x$ is no zero divisor
in $R/(x_1,...,x_s)$. We know that $\inn(x)$ is no zero
divisor in 
$\gr_I(R)/(\inn(x_1),...,\inn(x_s))\isom\gr_{\bar{I}}(\bar{R}).$ 
So by the first step, the assumption follows.
For the odd case we use the analogue argument.
So the claim is proven.

Now when (i) is true, it is clear that $\{\inn(x):\;x\in X\}$,
which is just the set $T$,
satisfies condition (ii) and by the claim, $X$ satisfies condition
(ii).\\

\textbf{Proof of (iv)$\impl$(i)}\\
Without restriction, we can suppose that $\sC$ is the category of commutative
$\QQ$-algebras.
For each $j\geq 0$,
we have to show that the j-th homogeneous component $(R/I[T])_j$ in the
$T$-grading of $R/I[T]$ maps isomorphically
to $I^j/I^{j+1}$.

We will already make use of the implication (ii)$\impl$(iii).
Set $S:=\QQ[T]$. 
We consider $R$ as $S$-algebra via the map $t(x)\mapsto x$. 
Obviously $T\sub S$ satisfies
condition (ii), so by (iii), the Koszul complex $K_S(T)$
is a DG-resolution of $\QQ$ over $S$.\\

We consider the exact sequence
$$0\nach(T)^j/(T)^{j+1}\nach S/(T)^{j+1}\nach S/(T)^j\nach 0$$
of graded $S$-modules.
$(T)^j/(T)^{j+1}$ is a free, graded g-finite $\QQ$-vectorspace
which is a $S$-module via the canonical map $S\nach\QQ$.
We write $\coprod_{i\in J}\QQ e_i$ for it.
Now $\coprod_{i\in J}K_S(T)e_i$ is a free resolution of 
$\coprod_{i\in J}\QQ e_i$ over $S$.
So we get 
$$\Tor_1^S((T)^j/(T)^{j+1},R)=
H^{-1}(\coprod_{i\in J}(K_S(T)e_i\ot_SR))=
\coprod_{i\in J}H^{-1}(K(X)e_j)=0.$$
By the property of left derived functors, there
is an exact sequence 
\begin{align*}
0\nach&\Tor_1^S(S/(T)^{j+1},R)\nach\Tor_1^S(S/(T)^{j},R)\nach\\
&(T)^j/(T)^{j+1}\ot_SR\nach S/(T)^{j+1}\ot_SR\nach S/(T)^j\ot_SR\nach 0
\end{align*}
By induction on j an the exactness of the first line, we
see that $\Tor_1^S(S/(T)^{j},R)=0$ for any $j\geq 0$.
The second line gives rise to a short exact sequence
$$0\nach(R/I[T])_j\nach R/I^j\nach R/I^{j+1}\nach 0,$$
which implies the desired isomorphism.\\

\textbf{Proof of (ii)$\impl$(iii)}\\
To prove this implication, we only have to show that for $p\leq 0$
the $p$-th row of
the double complex $K(X)$ is a DG-resolution in $\sM$ over $R^0$ of the p-th
component of $R/(X)$. For this we can suppose that $X$ is finite
with $g(x)\geq p$ for all $x\in X$.
Say $X=\{x_1,...,x_n\}$.
Then we have $K(X)=K(x)\ot...\ot K(x_n)$.\\

Each $K(X)^{(p,q)}$ is obviously a finite $R$-module, so with 
corollary~\ref{Apoc}.\ref{2506022} we 
only have to show that $K(X)^{(p,*)}$ is a resolution of $(R/(X))^p$ in 
the category of $R$-modules. To show this, we proceed by induction on m.\\
For m=1 we write $x$ instead of $x_1$ and $e$ instead of $e_1$.
Set  $m:=g(x)$.
Remark that if $m$ is even, then the total degree $m-1$ of $e$ is odd,
so in this case we have $e^2=0$. If $m$ is odd, then the total degree of $e$
is even, so $e^2\neq0$.
In the first case
$K(x)$ is just the complex

$$ \begin{array}{cccccccc}
   0 & \ldots & 0 & 0 & \ldots & 0 & 0& \ldots\\
   \downarrow & & \downarrow &\downarrow & &\downarrow &\downarrow &\\
   0 & \ldots & 0 & R^0e & \ldots & R^{m-1}e & R^m e & \ldots\\
   \downarrow & & \downarrow &\downarrow & &\downarrow &\downarrow &\\
   R^0 & \ldots & R^{m+1} & R^m & \ldots & R^{2m-1} & R^{2m} & \ldots \\
   \downarrow & & \downarrow &\downarrow & &\downarrow &\downarrow &\\
   R^0 & \ldots & R^{m+1} & R^m/R^0x & \ldots & R^{2m-1}/R^{m-1}x &
   R^{2m}/R^mx & \cdots
\end{array} $$

$s$ is injective since $x_1$ is no zero divisor in $R$, so
the sequence is exact.\\
In the second case
$K(x)$ is the complex

$$
\begin{array}{cccccccc}
   0 & \ldots & 0 & 0 & \ldots & 0 & R^0e^2 & \ldots\\
   \downarrow & & \downarrow &\downarrow & &\downarrow &\downarrow &\\
   0 & \ldots & 0 & R^0e & \ldots & R^{m-1}e & R^m e & \ldots\\
   \downarrow & & \downarrow &\downarrow & &\downarrow &\downarrow &\\
   R^0 & \ldots & R^{m+1} & R^m & \ldots & R^{2m-1} & R^{2m} & \ldots \\
   \downarrow & & \downarrow &\downarrow & &\downarrow &\downarrow &\\
   R^0 & \ldots & R^{m+1} & R^m/R^0x & \ldots & R^{2m-1}/R^{m-1}x &
   R^{2m}/R^mx & \cdots

 \end{array} 
$$

In $R^i$, for $i<m$ there is no element that annulates $x$,
so up to the row $m-1$, the situation is as above.
In the m-th row the kernel of $R^me\nach R^{2m}$ is just
$R^0xe$, so it coincides with the image of the map
$R^0e^2\nach R^me$.
Remark that here-fore we use that 2 is invertible in $R$.
Inductively we see that all rows are exact.
Here-fore we use that all naturals are invertible.\\

Now suppose that the statement is proven for m.
We set $L:=K(x_1,...,x_m)=R\lan e_1,...,e_m\ran$
and $K:=K(X)=K(x_1,\ldots, x_{m+1})$.
We write $x$ and $e$ instead of $x_m$ and $e_m$.
$K(x)$ is (as object of $\gr^2(\sM)(R)$) a direct product
$K_0\oplus K_{-1}\oplus K_{-2}\oplus\ldots$,
where in the case where $x$ is even, we have
$K_0=R$, $K_1=R[m,-1]$ and $K_s=0$ for $s<-1$ and
in the odd case we have
$K_s=R[sm,-s]$ for all $s\leq 0$.
The differential in $K(x)$ is given by the maps
$d_q:K^{p,q}_q\nach K^{p,q+1}_{q+1}$.
Remark that we have $K^{p,q}_q=R^{p+qm}$ and
$d_q:R^{p+qm}\nach R^{p+(q+1)m}$ is just multiplication by $x$.
Now $(K(x)\ot L)^{p,q}$ is the class
\begin{align*}
[\sum_{i+k=p,j+l=q}K_0^{i,j}\ot_{R^0}L^{k,l}\oplus
K_{-1}^{i,j}\ot_{R^0}L^{k,l}\oplus\ldots ]=\\
[\sum_{i+j=p}K_0^{i,0}\ot_{R^0}L^{k,q}\oplus
K_{-1}^{i,-1}\ot_{R^0}L^{k,q+1}\oplus\ldots   ]
\end{align*}
which is equal to 
$L^{p,q}\oplus L^{p-m,q+1}\oplus L^{p-2m,q+2}\oplus\ldots$
The differential of $K=K(x)\ot L$ is given by the scheme
\begin{equation*}
\xymatrix{
\sum_{i+j=p}K_0^{i,0}\ot_{R^0}L^{k,q}\ar[d]^{\tilde{\delta}}&\oplus&
\sum_{i+j=p}K_{-1}^{i,-1}\ot_{R^0}L^{k,q+1}\ar[dll]^x\ar[d]^{\tilde{\delta}}&
\oplus&\ldots\\  
\sum_{i+j=p}K_0^{i,0}\ot_{R^0}L^{k,q+1}&\oplus&
\sum_{i+j=p}K_{-1}^{i,-1}\ot_{R^0}L^{k,q+2}&\oplus&\ldots
}
\end{equation*}
But now in the even case the class
$[\sum_{i+k=p}K_s^{i,s}\ot_{R^0}L^{k,q-s}]$
is equal to $L[sm,-s]^{p,q}$ for $s=0,-1$ and $0$ for
$s<-1$.
In the odd case the class
$[\sum_{i+k=p}K_s^{i,s}\ot_{R^0}L^{k,q-s}]$
is equal to $L[sm,-s]^{p,q}$ for all $s\leq 0$.
So in the even case, the complex $K(X)^{p,*}$ is the total complex
of the double complex
\begin{equation*}
\xymatrix{
\vdots\ar[d] & \vdots\ar[d] & \vdots\ar[d] & \\
L^{p,-2}\ar[d]^\delta & L^{p-m,-2}\ar[d]^{\delta}\ar[l]^x & 0\ar[d]\ar[l] & 
\cdots\ar[l]\\
L^{p,-1}\ar[d]^{\delta} & L^{p-m,-1}\ar[d]^{\delta}\ar[l]^x & 0\ar[d]\ar[l] & 
\cdots\ar[l]\\
L^{p,0}& L^{p-m,0}\ar[l]^x & 0\ar[l] &\cdots\ar[l]
}
\end{equation*}

In the odd case $K(X)^{p,*}$ is the total complex of the double complex
\begin{equation*}
\xymatrix{
\vdots\ar[d] & \vdots\ar[d] & \vdots\ar[d] & \\
L^{p,-2}\ar[d]^\delta & L^{p-m,-2}\ar[d]^{\delta}\ar[l]^x & L^{p-2m,-2}
\ar[d]^\delta\ar[l]^x & 
\cdots\ar[l]\\
L^{p,-1}\ar[d]^{\delta} & L^{p-m,-1}\ar[d]^{\delta}\ar[l]^x & 
L^{p-2m,-1}\ar[d]^\delta\ar[l]^x & 
\cdots\ar[l]\\
L^{p,0}& L^{p-m,0}\ar[l]^x & L^{p-2m,0} \ar[l]^x &\cdots\ar[l]
}
\end{equation*}

The first double complex is a DDG-resolution in $\gr^2(\sM)(R^0)$
of the DG-module 
$$(R/(x_1,\ldots,x_n))^p\leftarrow (R/(x_1,\ldots,x_n))^p\leftarrow 0
\leftarrow\ldots,$$
where the left arrow stands for multiplication by $x$.
But this DG-module is a resolution 
of $(R/(x_1,\ldots,x_n,x))^p$ over $R^0$, since $g(x)$ is even. 
So $K(X)^{p,*}$ is a resolution of
$(R/(x_1,\ldots,x_n,x))^p$.
The second double complex is a DDG-resolution in $\gr^2(\sM)(R^0)$
of the DG-module 
$$(R/(x_1,\ldots,x_n))^p\leftarrow (R/(x_1,\ldots,x_n))^p
\leftarrow(R/(x_1,\ldots,x_n))^p \leftarrow\ldots,$$
where the arrows stand for multiplication by $x$.
But this DG-module is a resolution 
of \\$(R/(x_1,\ldots,x_n,x))^p$ over $R^0$, since $g(x)$ is odd. 
So $K(X)^{p,*}$ is a resolution of
$(R/(x_1,\ldots,x_n,x))^p$.
So for both cases the induction step is done.
\qed
\begin{bem}
The assumption $\QQ\sub\KK$ is used only to prove the implications 
(ii)$\impl$(iii) and (iv)$\impl$(i). The assumption that for each subset 
$Y\sub X$ we have $\cap_{n\geq 1} I^nR/(Y)=0$ is used only to
prove (i)$\impl$(ii). So if one wants to get rid of it, use
condition (ii) for the definition of regular sequences.
It can be stated in a slightly modified manner, which depends 
on the order of the elements of $X$, then. 
\end{bem}
\begin{defi}
Let $R$ be a DG-algebra in $\grcn$.
Let $(\alpha_i,u_i,g_i)_{i\in J}$ be a family in
$\sN\times U'$ and $X=\{x_i:\;i\in J\}$ a family of elements
with $x_i\in G_{u_i}(R_{\alpha_i}^{g_i})$
such that for $\beta,\beta'\sub\alpha$ the sets
$\{\rho_{\beta\alpha}(x_i):\;\alpha_i=\beta\}$ and
$\{\rho_{\beta'\alpha}(x_i):\;\alpha_i=\beta'\}$ are disjoint.
Suppose that
$$X_\alpha:=\cup_{\beta\sub\alpha}
\{\rho_{\beta\alpha}(x_i):\;\alpha_i=\beta\}$$ is a regular 
(resp. handy) ($s_\alpha$-)sequence 
in $R_\alpha$ for each $\alpha$.
Then $X$ is called a \textbf{regular (s-~)sequence} (resp.
handy (s-)sequence) in $R$.
\end{defi}

\begin{kor}
When $R=(R,s)$ is a DG-algebra in $\grcn$ and $X$ a handy s-sequence in $R$, 
then $K(X)$ is a DG-algebra in $\gr^2(\sC)^{\sN}$ and if
$X$ is regular, then $K(X)$ is a DG-resolution of $R/(X)$ over $R$.
\end{kor}

\textbf{Remarks:}
When $R$ carries the structure of a DG-algebra $(R,s)$,
one would like the Koszul complex to carry the structure of
a DDG-module. In general this is not the case.\\
When $X$ is an s-handy sequence then, since $I=(X)$ is s-stable, 
the algebra
$\gr_I(R)$ has the structure of a DG algebra in $\gr(\Alg)$, such that
each $I^n/I^{n+1}$ is a DG submodule of $\gr_I(R)$.
If for example $R$ is already a free DG-algebra in
$\gr(\QQ-\Alg)$ with a set $X$ of free algebra generators,
i.e. $R=R/I[X]$, then the differential of $\gr_I(R)=R$
differs in general from the differential $s$.
In this way we get a modified differential $\tilde{s}$ on $R$.
In a similar way we get a modified differential,
when $R$ is a free DG-algebra in $\gr(\sC)$, for any good pair
of categories $(\sC,\sM)$. This will play a role in section~\ref{Adt}.
In geometric language, going over from $s$ to $\tilde{s}$ is
a deformation to the normal cone.

\subsection{The universal module of differentials}\label{umod}

Let  $k\nach A$ be a morphism in $\gr(\sC)$.
Set $R:=A\ot_k^{\sC}A$ and denote the kernel of the multiplication map
$R\nach A$ in the category $\gr(\sM)(R)$ by $I$. In this
subsection we use some notations of subsection~\ref{barr}.\\

Attention: In general $A$ is the cokernel of the inclusion $I\incl R$
only in the category $R$-$\Mod$.
\begin{bem}\label{149021}
Suppose that $(\sC,\sM)$ satisfies (S1), all $A^i$ are finite
$A^0$-modules  and $I$ is a g-finite $R$-module,
Then we have
$A=R/I:=\Kokern(I\incl R)$. Here, of course, we mean the categorical cokernel.
\end{bem}

In \cite{BinKos}, (6.12) the \textbf{universal module $\Omega_{A/k}$
of $k$-differentials} is defined as the cokernel in $\gr(\sM)(A)$
of the map $b_2:A\ot_kA\ot_kA\nach A\ot_kA$, sending
$a\ot b\ot c$ to $ab\ot c-a\ot bc +(-1)^{bc}ac\ot b$ for
homogeneous elements $a,b,c\in A$.\\

Here we consider all tensor powers as $A$-modules with respect to the first 
factor.
Now when $(\sC,\sM)$ is $(\sC^{(0)},\sM^{(0)})$ or $(\sC^{(1)},\sM^{(1)})$, 
there is
a well-known canonical isomorphism $\Omega_{A/k}\isom I/I^2$
of DG-modules in $\gr(\sM)(A)$. Here we want to show that under
some weak assumptions this is true in good pairs of categories.
So suppose that the pair $(\sC,\sM)$ is good. For technical reasons, 
we also have to assume, that the marking $G$ is trivial. 
This is the case in all examples, when $\sM$ is non-graded.  
First of all, we have to ask if
we can consider the $R/I$-module $I/I^2$ as object of $\gr(\sM)(A)$.

\begin{bem}
Let $R=(R,s)$ be a DG-object in $\grc$ such that all $R^i$ are finite 
$R^0$-modules.
Consider an ideal $I\sub R$ which is
generated by a handy s-sequence $X=\{x_j:\;j\in J\}$ in $R$. Then $I$ is a 
DG-object of $\grm(R)$ and $I/I^2$ is isomorphic as $R/I$-module to a 
DG-object of $\grm(R/I)$.
\end{bem}
\bew
For each $x\in X$ we choose a free module generator $e(x)$ with
$g(e(x))=g(x)$ and we see that $I$ is the image\footnote{We remind that by
image we mean the kernel of the cokernel map.} of the map from
the free module $\coprod_{x\in X}Re(x)$ to $R$, defined by
$e(x)\mapsto x$. So $I$ is already an object of $\grm(R)$.
But since $I$ is $s$-stable by assumption, $I$ is a DG-module.

For each pair $i,j$ in $J$ with $i\leq j$ we choose a free module generator
$e_{ij}$ with $g(e_{ij})=g(x_i)+g(x_j)$. 
We get a homomorphism $\coprod_{i\leq j}Re_{ij}\nach R$
of modules in $\grm(R)$ by sending $e_{ij}$ to the product $x_ix_j$.
This homomorphism factorises through $I$, so there is a homomorphism
$\pi:\coprod_{i\leq j}Re_{ij}\nach I$ and 
there is an isomorphism of $R$-modules
$\Kokern\pi\isom I/I^2$. It it easy to see, that the differential $s$
induces a differential on $\Kokern\pi$, that makes it a DG-module in
$\grm(R)$.\\
Now $I/I^2$ is also an $R/I$-module
and in $R/I-\Mod$ the objects $\Kokern\pi$ and $\Kokern\pi\ot_R^{\grm}R/I$
are isomorphic. And the latter is an object of $\grm(R/I)$.
\qed

In the sequel, by $I/I^2$ in fact we mean $\Kokern\pi\ot_R^{\grm}R/I$.\\

\begin{prop}
Let $k\nach A$ be a homomorphism of DG-objects in $\grc$. Suppose that
all $A^i$ are finite $A^0$-modules and that 
$I:=\Kern(\mu:A\ot_k^{\sC}A\nach A)$ is generated by an $s$-handy 
sequence $X$ in $ R:=A\ot_kA\nach A$. 
Here $s$ denotes the 
differential of $R$, induced by the differential of $A$.
Then by $\bar{a}\mapsto [a]$ we get an isomorphism
$I/I^2\nach\Omega_{A/k}$ in $\grm(R)$, whose inverse is given by
$[\alpha]\mapsto\overline{\alpha-\iota_1\mu(\alpha)}$.
Here $\overline{a}$ denotes the class in $I/I^2$ represented by $a$
and $[\alpha]$ denotes the class in $\Omega_{A/k}$ represented by
$\alpha$.
\end{prop}
\bew
First we have to show that the map $I/I^2\nach\Omega_{A/k}$,
$\bar{a}\mapsto [a]$ is well defined.\\
There is a homomorphism $\eta:I\ot_RI\nach I$ in $\grm(R)$
with $a\ot b\mapsto ab$. Consider the homomorphism
$\xi:I\nach \Omega_{A/k}$, $a\mapsto [a]$. For the well-definedness
it is enough to prove that $\xi\circ\eta=0$.
Since the barcomplex $C^{\ba}_{\bullet}(A/k)$ is acyclic, we see that
$b'$ gives rise to an epimorphism $A^{\ot^3}\nach I$.
Hence it is enough to show that the map
\begin{align*}
A^{\ot^6}&\nach \Omega_{A/k}\\
a\ot b\ot c\ot d\ot e\ot f&\mapsto [(ab\ot c-a\ot bc)(de\ot f-d\ot ef)]
\end{align*}
is zero.
But the argument in the brace on the right hand-side is just 
the image of
\begin{align*}
(-1)^{cd+ce+db}adbe&(cf\ot1\ot1-1\ot c\ot f)-\\
(-1)^{bd+cd}adb&(cef\ot1\ot1-1\ot c\ot ef)-\\
(-1)^{bd+be+cd+ce}ade&(bcf\ot1\ot1-1\ot bc\ot f)+\\
(-1)^{bd+cd}ad&(bcef\ot 1\ot1-1\ot bc\ot ef)
\end{align*} 
by the map $b_2$.

Secondly we have to show that the map
$\Omega_{A/k}\nach I/I^2$, $[\alpha]\mapsto\overline{\alpha-s_1\mu(\alpha)}$
is well defined. But there is a derivation
\begin{align*}
\delta:A&\nach I/I^2\\
a&\mapsto \overline{1\ot a-a\ot1}
\end{align*}
So by the universal property of $\Omega_{A/k}$,
see the proof of \cite{BinKos}
 lemma (6.13), there is a map $\Omega_{A/k}\nach I/I^2$
sending a class $[a\ot b]$ to $a\delta(b)=
a\cdot \overline{1\ot b-b\ot1}$ and we see that this map is just the
map we want.

To see that the both given maps are
inverse to each other, we remark that
elements of the form $a\ot 1$ in $A\ot A$ are in the image
of $b_2$, so they represent the zero class.
\qed

$I/I^2$ has the structure of an $A$-module
in $\grm(A)$. The multiplication $A\times I/I^2\nach I/I^2$
is inherited by the multiplication
$a\cdot\alpha=\iota_1(a)\cdot\alpha$ on $A\ot A$.
But on $A\ot A$ there is also a left multiplication
$\alpha\cdot a:=\alpha\cdot\iota_2(a)$.
Remark that the left- and right multiplication induced on
$I/I^2$ make $I/I^2$ an antisymmetrical $A$-bimodule.\\ 

Now, let $R=(R,s)$ be a DG-object of $\grc$ and suppose that all $R^i$ are
finite $R^0$-modules. Let $I\sub R$ be an ideal which is
generated by a regular s-sequence $X\sub R$. Say $X=\{x_i:\; i\in J\}$.
s defines a differential $\delta$ on $R/I$, that we denote again by $s$.
Consider the free module $\coprod_{i\in J}R/Ie_i$, where
the $e_i$ are free module generators of degree $g(x_i)$.
\begin{bem} 
We can make $\coprod_{i\in J}R/Ie_i$ a DG-module, by
defining a differential $t$ in the following sense:
For $i\in J$, there is a finite family of elements $x_{ij}\in X$
and $a_{ij}\in R$ such that $s(x_i)=\sum_ja_{ij}x_{ij}$.
Now we set $\delta(a):=s(a)$ for elements $a\in R/I$ and
$\delta(e_i):=\sum_j\overline{a_{ij}}e_{ij}$.
\end{bem}
\bew
To show that this defines a differential on $\coprod_{i\in J}R/Ie_i$,
we only have to show, that $\delta^2(e_i)=0$ for $i\in J$.
But since $s$ is a differential on $R$, we have
$$0=s^2(x_i)=s(\sum_ja_{ij}x_j)=
\sum_{j,k}(-1)^{a_{ij}}a_{ij}a_{jk}x_k+\sum_js(a_{ij})x_j.$$
We can reorganize the coefficients and get a sum
$\sum_{k=1}^mb_kx_{i_k}=0$ where the $x_{i_k}$ are pairwise different.
Remark that $\sum_{k=1}^m\overline{b_k}e_{i_k}=0$ is just $\delta^2(e_i)$.
To show that this sum is zero, we have to show that each $b_k$ belongs to
$I$. But assume that one $b_k$, say $b_m$ does not belong to $I$, then
$\overline{b_m}$ is a nonzero annulator of $x_m$ in $R/(x_1,\ldots,x_{m-1})$
and it doesn't belong to $Rx_m$. This contradicts the hypothesis
that $X$ is regular.
\qed

In the algebraic case, the following proposition is an immediate
consequence of condition (i) in definition and theorem~\ref{Apoc}.\ref{regS}.

\begin{prop} In this situation, the assignment
\begin{align*}
\coprod_{x\in X}R/I e(x)  \nach I/I^2,\quad
e(x)\mapsto \bar{x}
\end{align*}
gives rise to an isomorphism of DG-objects in $\grm(R/I)$.
\end{prop}

\bew
It is clear that the map commutes with the
differentials.\\
Obviously the map is well defined and surjective.
With axiom (S2), we only have to show that the map is
injective.
So let $\sum_{i=1}^m e_i\bar{a}_i$ be an element of the kernel of this map.
Then we have $\sum\overline{a_ix_i}=0$, i.e.
$\sum a_ix_i\in I^2$.
We must show that all $a_i$ are elements of $I$.
Let $Y$ be a finite subset of $X$ such that $\sum a_ix_i$ is
a sum $\sum_{y,y'\in Y}a(yy')yy'$ with $a(yy')\in R$.
Now as in the well-known nongraded case, when we assume that
one $a_i$, say $a_m$ is not in $I$, we can deduce that
$a_m$ is a zero divisor in $R/J$, where $J\sub R$ is the ideal
generated by $Y\setminus x_m$. This leads to a contradiction!
\qed

%\begin{defi}
%$A$ is called \textbf{very smooth} over $k$
%if the kernel of the
%multiplication map $R\nach A$ is generated by a regular sequence $X$ in $R$.
%\end{defi}

The condition on $A$ in the following corollary
is something like a smoothness condition.

\begin{kor}\label{Omega}
Suppose that all $A_i$ are finite $A^0$-modules.
If the kernel of the
multiplication map $R:=A\ot_kA \nach A$ is generated by a regular 
s-sequence $X$ in $R$
then there is a natural isomorphism of DG-modules
in $\grm(A)$
\begin{equation*}
\Omega_{A/k}\nach\coprod_{x\in X}A e(x).
\end{equation*}
Here $X$ denotes the regular s-sequence in $R$ and to $x\in X$ we
have associated a free module generator $e(x)$ with $g(e(x))=g(x)$.
The differential on the left 
is given by the rule $e(x)\mapsto\sum\bar{a}_ye(x_y)$,
where for $x\in X$ the family $a_y$ is chosen in a way, such that
$s(x)=\sum a_yy$ and $\bar{a}$ denotes the residue class in
$R/(X)\isom A$ of an element $a\in R$.
\end{kor}

From this statement, we can deduce the corresponding simplicial
statement:
\begin{prop}\label{Omegan}
Suppose that $A$ is a DG-algebra in $\grcn$ over $k$ and
set $R:=A\ot_kA$. Denote the differential on $R$, induced by 
the differential of $A$ by $s$.
Suppose that the kernel of the map $R\nach A$ is generated by
a regular $s$-sequence $X$ in $R$ in the sense that for each
$\alpha\in\sN$, the kernel of $R_\alpha\nach A_\alpha$ 
is generated by 
$X_\alpha:=\cup_{\beta\sub\alpha}\{\rho_{\beta\alpha}(x_i):\;\alpha_i=\beta\}$.
Then there is an isomorphism in the category of 
DG-modules in $\grmn$
$$\Omega_{A/k}\nach \coprod_{x\in X}Ae(x).$$
\end{prop}

\vspace{-1cm}
\section{Hochschild complexes and Hochschild cohomology}\label{hchc}

\subsection{The algebraic context}\label{tac}

In this subsection we want to generalize the algebraic
definition of Hochschild complexes in such a way that it is also
useful in analytical contexts. Here let $\sC$ be the
category of commutative rings and $\sM$ the category of
modules over objects of $\sC$.\\

Let $k\nach a$ be a homomorphism in $\sC$.
Let $A$ be a resolvent of $a$ over $k$.
Set $R:=A\ot_k A$.
As in subsection~\ref{barr}, 
we denote the cyclic bar complex and the bar complex of a
DG-algebra over $k$
by $C^{\cycl}_\cdot(\;)$ and $C^{\ba}(\;)$.

The double complex $C^{\cycl}(A)$ is a DDG-resolution of
$C^{\cycl}(a)$. Hence the free $R$-algebra
$\tot(C^{\cycl}(A))$ is quasi-isomorphic over
$R$ to $C^{\cycl}(a)$. Since $C^{\ba}(A)\ot_RA\isom C^{\cycl}(A)$,
there are quasi-isomorphisms:
$$C^{\cycl}(a)\approx\tot(C^{\cycl}(A))\approx\tot(C^{\ba}(A))\ot_Ra.$$

Now since in the algebraic case,
$\tot(C^{\ba}(A))$ is a resolvent of $A$ over $R$, and two
such resolvents are homotopy-equivalent, we get the following
result:
\begin{prop}\label{hochs}
Let $S$ be a g-finite resolvent of $A$ over $R$.
Then there is a quasi-isomorphism
$$S\ot_Ra\nach C^{\cycl}(a)$$  over $a$.
\end{prop}

Remember that on the right hand-side we have the classical Hochschild
complex. This shall justify the definition in the next subsection.
Remark that proposition~\ref{hchc}.\ref{hochs} keeps true in the simplicial
context.

\subsection{The noetherian context}\label{secta}

Fix a good pair of categories $(\sC,\sM)$ with marking $(F,G)$,
where $G$ is the trivial marking of $\sM$. Suppose that the axioms
(N) and (F2) are satisfied.\\

Let $k\nach a$ be a finite morphism of $\sN$-objects in $\sC$, i.e.
$a$ is a quotient of a free $k$-algebra $b$ in $\sC^{\sN}$, such that 
for each $\alpha\in\sN$ the algebra $a_\alpha$ is a 
a free finite $k_\alpha$-algebra.
Then, with \cite{BinKos}, proposition (8.8), there exists a g-finite
resolvent of $a$ over $k$. Fix such a resolvent $A$.
Set $R:=A\ot^{\grc}_k A$ and consider $A$ as algebra over $R$ by
the multiplication map $\mu:R\nach A$.
Let $S$ be a free g-finite
resolvent\footnote{Again with loc. cit., such a resolvent exists.
We can even construct it in such a way that $S^0=R^0$.} of $A$ over $R$.

\begin{defi}
We define the \textbf{simplicial Hochschild complex} $\Hk_{\ast}(a/k)$
of $a$ over $k$ to
be the object represented by the complex $S\ot_R a$ in the
homotopy category $K^{-}(\sM^{\sN}(a))$.
\end{defi}
\vspace{-0.7cm}

\begin{prop}\label{hhoinv}
$\Hk_{\ast}(a/k)$ is a well defined object in $K^{-}(\sM^{\sN}(a))$.
\end{prop}
\bew
For $i=1,2$ let $A_i$ be a g-finite resolvent of $a$ over $k$,
$R_i:=A_i\ot_kA_i$ and let $S_i$ be a g-finite resolvent of $A_i$ over
$R_i$. We have to show the existence of a homotopy equivalence
$S_1\ot_{R_1}a\simeq S_2\ot_{R_2}a$ over $a$.\\
First remark that $R_i$ is a g-finite resolvent of $a\ot_k a$ over $k$.
Hence by proposition~\ref{Apoc}.\ref{rehoeq}, 
there is a homomorphism $R_1\nach R_2$
in $\grcn$ which is a homotopy equivalence over $k$.
Hence we get a quasi-isomorphism
$$S_1=S_1\ot_{R_1}R_1\nach S_1\ot_{R_1}R_2$$ over $R_1$.
$S'_1:=S_1\ot_{R_1}R_2$ is a free algebra over $R_2$ and
as $S_2$ it is a resolution of $a$ over $R_2$.
Hence by proposition~\ref{Apoc}.\ref{rehoeq}, there is a 
homomorphism $S'_1\nach S_2$
in $\grcn$, which is a homotopy equivalence over $R_2$.
We can tensorise both sides over $R_2$ with $a$ and still get
a homotopy equivalence
$S_1\ot_{R_1}a\nach S_2\ot_{R_2}a$.
\qed

We consider $a$ as object of $\gr(\sM)^{\sN}(A)$ via the surjection
$\alpha:A\nach a$. In this sense the Hochschild complex can be seen as
DG-module over $A$ in $\gr(\sM)^{\sN}$.

\begin{bem}\label{jako}
The map $\tilde{\alpha}:S\ot_R A\nach S\ot_R a$ induced by 
$\alpha$ is a quasi-isomorphism over $A$. 
\end{bem}
\begin{defi}
Let $M$ be an object of $\sM^{\sN}$ over $a$.
We define the \textbf{Hochschild cochain complex} of $a$ over $k$ 
with values in $M$ to be the complex
$$\Hom_a^{\sN}(\Hk_{\ast}(a/k),M),$$ 
with the differential induced by the differential
of $\Hk_{\ast}(a/k)$.
We define the \textbf{Hochschild cohomology} $\HH(a/k,M)$ of $a$ over $k$ 
with values in $M$ to be 
the cohomology of the Hochschild cochain complex.
\end{defi} 
\begin{prop}
The Hochschild cochain complex is well defined up to homotopy equivalence.
\end{prop}
\bew
This is a consequence of proposition~\ref{hchc}.\ref{hhoinv} and
\cite{BinKos}, chapter I, lemma (3.7).
\qed
\begin{lemma}\label{homme}
The Hochschild cohomology of $a$ over $k$ with values in $M$ is equal to the
cohomology of the complex $\Hom_{A}^{\sN}(S\ot_{R} A,M)$, where $M$ is 
considered as object over $A$ via the map $A\nach a$.
\end{lemma}
\bew
With \cite{BinKos} lemma 6.4, we have
\begin{equation}\label{glei} 
\Hom_a(a\ot_R S,M)=\Hom_R(S,M)=\Hom_A(A\ot_R S,M).
\end{equation}
\qed

\section{A decomposition theorem for Hochschild cohomology}\label{Adt}

Let $A=k\langle T\rangle$ be a g-finite free DG-object of $\grc$ over $k$,
with differential $d$. We think of $A$ as the resolvent of a
$k$-algebra $a$ in $\grc$.
Further we consider the free algebra $R:=A\ot_kA$ in $\grc_k$ with
its differential $s=d\ot 1+1\ot d$. 
Then the multiplication map $\mu:R\nach A$ is a homomorphism of
DG-objects in $\grc$.\\
For the construction of the Hochschild complexes of $a$ over $k$, we
are interested in resolvents $S$ of $A$ as DG-object in $\grc(R)$.\\

First, we will not construct such an $S$
but something very similar. 
For this we modify the differential $s$ on $R$
as sketched in subsection~\ref{hans}, which geometrically corresponds
to a deformation to the normal cone.
We call the modified differential $\tilde{s}$.
This modification doesn't touch the structure of
$A$ as DG-algebra over $R$, since $s(r)-\tilde{s}(r)$
will be in the kernel of $\mu$ for each $r\in R$.

\newcommand{\ts}{\tilde{s}}
Then 
we construct a resolvent $S$ of $A$ over $\tilde{R}:=(R,\ts)$
with the help of a Koszul complex. As we will see in
subsection~\ref{auflS}, this is enough 
for the construction of the Hochschild complex.\\

If $B$ is an object of $\grc$ and $R:=B\langle T\rangle$
is a free algebra over $B$ in $\grc$ with a g-finite set
$T$ of
free generators $t$ with $t\in F_{\tau(t)}(R^{g(t)})$.
Then $R\ot_BR$ is a free algebra over $B$ with
two free algebra generators $t_1=t\ot 1$ and $t_2=1\ot t$ for each
$t\in T$. For $t\in T$ set $t^+:=\frac{1}{2}(t_1+t_2)$ and
$t^-:=\frac{1}{2}(t_1-t_2)$. Let $T^+$ be the set of all $t^+$
and $T^-$ be the set of all $t^-$.

We say that the marking $F$ on $\sC$ is \textbf{balanced}
if for each $\tau\in\T$ and each $A$ in $\sC$ and each $t\in
F_\tau(A)$ we have $-t\in F_\tau(A)$. We say that the marking
$F$ is \textbf{convex} if for each $\tau\in\T$, each $A$ in $\sC$,
each $t_1,t_2\in F_\tau(A)$ and each $a,b\in\KK$ with $a+b=1$
we have $at_1+bt_2\in F_\tau(A)$. 

\begin{bem}\label{axiom}
Suppose that the marking $\gr_G(F)$ on $\grc$ is balanced
and convex, then we have
$R\ot_B R\isom B\langle T^{+}\cup T^{-}\rangle$.
More precisely this means that there is a free algebra over $B$ and
an isomorphism from this free algebra in $R\ot_BR$, sending
the free generators on the elements $t^+$ and $t^-$.
\end{bem}
\begin{beisp}
\begin{liste}
\item
The trivial marking on $\sC$ is balanced and convex,
so when $\sC$ is the category $\sC^{(0)}$ of (noetherian) rings
and $\sM$ is the category $\sM^{(0)}$
of modules over $\sC^{(0)}$, then remark~\ref{Adt}.\ref{axiom}
is true.
\item
When $\sC$ is the category $\sC^{(1)}$ of (local) analytic algebras and
$\sM$ the category $\sM^{(1)}$ of DFN-modules over $\sC^{(1)}$, 
then the marking
$F$ on $\sC$ (see example~\ref{Apoc}.\ref{makke})  is balanced and convex.
\end{liste}
\end{beisp}
\bew
The first example is trivial. For the second example, we show that
if a free generator $t$ is in $F_{\tau}(R)$, then
$T^+$ and $t^-$ are in $F_{\tau}(R\ot_BR)$:
Here $\tau$ stands for a positive real number
and $F_\tau(R)$ is the set of all $r\in R$, such that for each
character $\xi\in\sX(R)$, we have $|\xi(r)|\leq\tau$.
Now $t_1=\iota_1(t)=t\ot 1$ and $t_2=\iota_2(t)=1\ot t$ belong to
$F_\tau(R\ot_BR)$, so for each character $\xi\in\sX(R\ot_BR)$, we
have $|\xi(t_1)|\leq\tau$ and $|\xi(t_2)|\leq\tau$. Hence
$|\xi(t^+)|=|\frac{1}{2}(\xi(t_1)+\xi(t_2))|\leq\tau$ and
$|\xi(t^-)|=|\frac{1}{2}(\xi(t_1)-\xi(t_2))|\leq\tau$.
The case of local analytic algebras is clear, since
maximal ideals are additively closed.
\qed

\subsection{Deformation to the normal cone}

For the rest of this section suppose that the marking
$F$ on $\sC$ is balanced and convex and that the marking $G$
on $\sM$ is trivial. Further suppose that (N) and (F2) are satisfied.

Then we have $R=k\lan T^+\cup T^-\ran$. Since $T^-$ is $s$-stable
and g-finite, it is a regular $s$-sequence.
(In the algebraic case we have $R=\gr_{(T^-)}(R)$, so
as we have already mentioned in section~\ref{hans}, 
there is a deformation $\tilde{s}$ of $s$, that
respects the submodules $(T^-)^j/(T^-)^{j+1}$. 
This is what we mean by ``deformation to the normal cone''.
Here we do a similar construction for the general case.)
Each element $r$ of $R$ has a unique decomposition 
$r=\check{r}+\dot{r}+\hat{r}$ with
$\check{r}\in\check{R}:=k\lan T^+\ran$,
$\dot{r}\in\dot{R}:=\sum_{t\in T^-}t k\lan T^+\ran$ and
$\hat{r}\in \hat{R}:=\sum_{t,t'\in T^-}tt'k\lan T^+\ran$.
Now we define a $R^0$-derivation $\tilde{s}$ on $R$
setting $\tilde{s}(t):=(s(t))^\vee$ for $t\in T^+$ and
$\tilde{s}(t):=(s(t))^\cdot$ for $t\in T^-$.
The philosophy of this modification is that 
roughly speaking $\tilde{s}$ preserves
the $T^-$-degree of homogeneous elements in $R$. More precisely
we have $\tilde{s}(\check{R})\sub\check{R}$,
$\tilde{s}(\dot{R})\sub\dot{R}$ and $\tilde{s}(\hat{R})\sub\hat{R}$ in contrast
to $s(\dot{R})\sub\dot{R}\amalg\hat{R}$ and 
$s(\hat{R})\sub\hat{R}$.
\begin{prop}
$\tilde{s}$ is a differential, i.e. $\tilde{s}^2=0$. 
\end{prop}
\bew
First remark that for $a\in k\lan T^+\ran$ we have 
$$(s(a))^\vee=\tilde{s}(a).$$ To prove this we can suppose that
$a$ is of the form $a_0t_1\ppp t_n$ with $a_0\in k\lan T^{+,0}\ran$ and
$t_i\in T^{+,<0}$. In this case it is easy to see.\\
Now suppose that $t$ is in $T^+$. Then $s(t)=\tilde{s}(t)+\text{rest}$,
where $\text{rest}$ is in $\dot{R}\amalg\hat{R}$.
So $s^2(t)=\tilde{s}^2(t)+\text{rest}'$
where $\text{rest}'$ is in $\dot{R}\amalg\hat{R}$.
Since $\tilde{s}^2(t)$ is in $\check{R}$ and $s^2(t)=0$,
we get $\tilde{s}^2(t)=0$. Similarly, we see that $\tilde{s}^2(t)=0$
for $t\in T^-$, which proves the proposition.
\qed

We write $X$ for the regular $\tilde{s}$-sequence $T^-$.
Next we will see that for $\tilde{R}=(R,\tilde{s})$, 
the Koszul complex $(K(X),v)$ has the structure $(K(X),h,v)$
of a DDG-algebra, so
its total complex is 
a resolution of $A=\tilde{R}/(X)$ over $\tilde{R}$:
Again we denote by $E$ the set of free algebra generators,
containing for each $x_i\in X$ an element $e_i$ of bidegree
$(g(x_i),-1)$. Here $\tilde{s}(x_i)$ is a sum of the form
$\sum a_j x_j$,
where no $a_j$ belongs to the ideal $(X)$. So in fact all $a_j$ belong
to $B=k\langle T^{+}\rangle$.
Now there is exactly one choice for the element $h_i$, which shall
be the image of $e_i$ by the horizontal differential  $h$
of $ R\langle E\rangle $. The choice is $h_i=\sum a_je_j$.
Now we have 
$0=\tilde{s}^2(x_i)=\tilde{s}(\sum a_jx_j)=\sum_{j,k}a_ja_{jk}x_{jk}$   
and the coefficients $a_j a_{jk}$ belong to $k\langle T^{+}\rangle$.
So we have 
$h^2(e_i)=h(\sum a_je_j)=\sum_{j,k}a_ja_{jk}e_{jk}=0$.
I.e. the hypothesis of remark~\ref{Apoc}.\ref{DDGstr} is satisfied.
So the Koszul complex $(K(X),v)$, equipped with 
the horizontal differential $h$,
is a DDG-resolution in $\gr^2(\sC)$ of $A=R/(X)$ over 
$\tilde{R}=(R,\tilde{s})$. And the total complex $\tilde{S}$
of $K(X)$ is a free algebra resolution of $A$ over $\tilde{R}$.\\

The reason for this construction is that there is a nice
description of the tensor product 
$$\tilde{S}\ot_{\tilde{R}}A=\tot(K(X))\ot_{\tilde{R}}A=
 \tot(K(X)\ot_{\tilde{R}}A).$$ Namely:
\begin{bem}\label{til}
The double complex $K(X)\ot_{\tilde{R}}A$
is just $\wedge^\cdot_A\Omega_{A/k}$.
\end{bem}
\bew
This follows imidiately by proposition~\ref{Apoc}.\ref{Omegan} 
and definition~\ref{Apoc}.\ref{dach}.
\qed

The going over from $s$ to $\tilde{s}$ is natural, i.e.
when $(R,s)=(r_\alpha,s_\alpha)_{\alpha\in\sN}$ is a DG-algebra
in $\grcn$, then
$\tilde{R}:=(R_\alpha,\tilde{s}_\alpha)_{\alpha\in\sN}$
is again a DG-algebra in $\grcn$ and remark~\ref{Adt}.\ref{til}
keeps true in the simplicial case.

\subsection{Construction of the resolution $S$}\label{auflS}

To get a resolution $S$ of $A$ over $R$, which is good for our purpose,
we have to work harder. The strategy is to construct again a DDG-resolution
$K$ in $\gr^2(\sC)^{\sN}$, that is free over the double graded object $K(X)$
and such that the projection $K\nach K(X)$ is a morphism of DG-
objects (of course it will not be a morphism of DDG-objects)
and such that the induced map $K\ot_RA\nach K(X)\ot_{\tilde{R}}A$
is a morphism of DDG-objects in $\gr^2(\sC)^{\sN}$, 
that induces an isomorphism
on vertical homology. When we have managed to realize this,
with $S:=\tot(K)$ we get a resolution of $A$ over $R$ and a 
quasi-isomorphism 
\begin{equation}\label{heck}
S\ot_RA\nach \tilde{S}\ot_{\tilde{R}}A
\end{equation}
over $A$.
Since both complexes are free objects over $A$, the quasi-isomorphism is even
a homotopy-equivalence.\\ 

First we explain heuristically the construction of $K$:
We take the Koszul-complex $K(X)$ with its vertical differential $v$.
The problem is to define a horizontal differential $h$ on it, since here
in general there are no good candidates for the values
$h(e)$ of $h$ on the free generators $e\in E$. 
When we have $s(x)=\sum a_{xy}y$, to get commutative diagrams,
$h(e(x))$ must be something like $\sum a_{xy}e(y)$. The problem is
that $h(\sum a_{xy}e(y))$ won't be zero.
So, inductively, for $e$ we add free algebra-
generators $f$ of bidegree $(g(x)+1,-1)$ with $v(f)=0$,
in such a way that we can find candidates for $h(e)$ in
$K(X)^{g(x)+1,-1}+\sum R^0f$. 
When this is done for all $e\in E$, 
we get a DDG-structure on the extension $K(X)\lan F\ran=R\lan E\cup F\ran$.
Now this extension is not any more a resolution. To get 
a resolution again, we apply the construction of remark~\ref{Apoc}.\ref{fres}\\

Now we begin with the construction of $K$. First we place ourselves
in the affine situation.\\
For each $x\in X$ we fix a finite family $a_{xy};\;y\in Y\sub X$
such that $s(x)=\sum a_{xy}y$.\\
Set $K(X):=R\lan E\ran$ as double graded algebra.

\begin{prop}\label{beha1}
There is a g-finite family $F=\cup_{p\leq 0}F^p$
of free algebra generators with $g(f)=(p,-1)$ for $F\in F^p$
and a DDG-algebra
structure $(L,h,v)$ on $L:=R\lan E\cup F\ran$ such that
\begin{liste}
\item
$v(e(x))=x$
\item\label{beha12}
$v(f)=0$
\item
$h(f)$ is in the ideal $(X\cup F)$, generated by $X$ and $F$.
\item
$h(e(x))=\sum_ya_{xy}e(y)+\gamma(x)$ with a $\gamma(x)\in (F)$.\label{beha14}
\end{liste}
\end{prop}
\bew
We construct a sequence $L_k=(L_k,h_k,v_k)$ of free DDG-algebras
over $R$ as well as a family $\{\gamma(x):\;x\in X\wedge g(x)\geq k-1 \}$
with $\gamma(x)\in L_k^{k,-1}$ for $x\in X^{k-1}$ such that the following 
conditions hold:
\begin{enumerate}
\item[(a)]
$L_0=R\lan E^0\ran$.
\item[(b)]
$L_{k-1}=L_k\lan E^{k-1}\cup F^{k-1}\ran$ is a free DDG-algebra over $L_k$,
where $F^{k-1}$ is a finite set of algebra generators of bidegree
$(k-1,-1)$ and we have the rules $v_{k-1}(e(x))=x$ (as in the Koszul-construction) and $v_{k-1}(f)=0$.
\item[(c)]
$h_{k-1}$ maps the submodule $\amalg_{f\in F^{k-1}}R^0f$ of $L_{k-1}^{k-1,-1}$
surjectively onto $\Kern(h_k^{k,-1})\cap\Kern(v_k^{k,-1})\sub L_k^{k,-1}$.
\item[(d)]
For all $i\geq k-1$ the sequence $L_{k-1}^{i,-1}\nach R^i\nach A^i$ in
$\sM(R^0)$ is exact.
\item[(e)]
For $x\in X^{k-1}$ we have $\gamma(x)\in\amalg_{f\in F^k}R^0f$ and
$h_k^{k,-1}(\gamma(x))=h_k^{k,-1}(\sum_y a_{xy}e(y))$
\end{enumerate}

For $L_0=R\lan E^0\ran$ we have already seen that by setting
$v(e(x)):=x$ we get an exact sequence $L_0^{0,-1}\nach R^0\nach A^0$.
Now suppose that $L_k$ and $\{\gamma(x):\;g(x)\geq k\}$ is already
constructed.
We choose finitely many free algebra generators $f$ of bidegree
$(k-1,-1)$ such that there exists
an epimorphism
$$\pi:\amalg R^0f\nach\Kern(h_k^{k,-1})\cap\Kern(v_k^{k,-1})\sub
(X\cup F).$$
To explain the inclusion: The vertical differential on the subalgebra
$K(X)$ is exact. So a homogeneous element of $K(X)\lan F\ran$, which is in
the kernel of $v$, is a sum of an element in the image of $v$ and an element
in the ideal (F).\\
We set $h_{k-1}(f):=\pi(f)$ and $v_{k-1}(f):=0$.
For $x\in X^{k-1}$, to see that there exists a good candidate for
$\gamma(x)$, we must show that $h_k(\sum_y a_{xy}e(y))$ belongs to
$\Kern(h_k^{k,-1})\cap\Kern(v_k^{k,-1})$. But we have
$h_k(\sum_y a_{xy}e(y))=\sum_ys(a_{xy})e(y)+
\sum_ya_{xy}[\sum_za_{yz}e(z)-\gamma(y)]=
(\sum_ys(a_{xy})e(y)+\sum_{y,z}a_{xy}a_{yz}e(z))-\sum_ya_{xy}\gamma(y)$.
The first term maps vertically to $s^2(x)$ which is zero, the second
factor maps obviously vertically to zero.\\

At least we can set $h_{k-1}^{k-1,-1}(e(x)):=\sum a_{xy}e(y)-\gamma(x)$
and $v_k^{k-1,-1}(e(x)):=x$. This gives the desired $L_{k-1}$.
Then we set $L:=\lim L_k$.
\qed

\textbf{Properties of $L$:}
\begin{liste}
\item
$L^{0,*}=K(X)^{0,*}$, hence this is a resolution of $A^0$ over $R^0$.
\item
$L^{0,p}=R^p$ for all $p\leq 0$.
\item
The sequence $L^{-1,p}\nach L^{0,p}\nach A^p\nach 0$ is exact for all 
$p\leq 0$.
\item
The inclusion $K(X)\incl L$ and the projection $L\nach K(X)$ are
homomorphisms of DG-algebras over $(R,0)$, so in the category of
DG-modules in $\gr(\sM)$ there is a decomposition $L=K(X)\amalg L'$.
The (vertical) homology of $L$ is contained in $L'$.
\end{liste}

\begin{prop}\label{beha2}
There is a g-finite
family $G=\cup_{p\leq =,q\leq -2}G^{p,q}$ of free algebra generators
with $g(g)=(p,q)$ for $g\in G^{p,q}$ and extensions of $h$ and $v$
on $K:=L\lan G\ran$, such that
\begin{liste}
\item
The i-th row of $K$ is a $R^0$-module resolution of $A^i$.
\item
$v(g)\in (F\cup G)$\label{beha22}
\item
$h(g)\in (G)$
\end{liste}
\end{prop}
\bew
We can construct the free DDG-resolution $K$ of $A$ over $R$
with the method of remark~\ref{Apoc}.\ref{fres}.
\qed

Comparing the values of $h$ on the free generators $e$ with
its values by the differential $\tilde{h}$ of the Koszul-complex
over $R$ with the modified differential, we see that
$\tilde{h}(e)-h(e)\in \amalg R^0g+\sum_{x\in X}xK$.\\
 
\textbf{Consequence:}
Consider the projection
$\pi:K=R\lan E\cup F\cup G\ran \nach R\lan E\ran=K(X)$
(a priori only as map of algebras in $\gr^2(\sC)$).
With proposition~\ref{Adt}.\ref{beha1} (ii) and
proposition~\ref{Adt}.\ref{beha2} (ii), we see that $\pi$
respects the vertical differential. By the construction of $\tilde{s}$
and proposition~\ref{Adt}.\ref{beha1} (iv), we see that
$$\pi\ot 1:K\ot_RA\nach K(X)\ot_{\tilde{R}}A$$
is a homomorphism of DDG-algebras in $\gr^2(\sC)$ over $A$.\\

Now we can prove the (affine case of the) crucial result of this chapter.
It says that to construct the Hochschild complex
it is enough to work with a resolvent of $A$ over $\tilde{R}$. 
\begin{satz}\label{hohoho}
With $S:=\tot(K)$ and $\tilde{S}:=\tot(K(X))$, there is a 
homotopy-equivalence
$$S\ot_RA\nach\tilde{S}\ot_{\tilde{R}}A$$
over $A$.
\end{satz}

\bew
First we have seen that the projection $\pi:K\nach K(X)$ is
a homomorphism of DG-Algebras in $\gr^2(\sC)$ over $(R,0)$.
Since both double complexes are free resolutions, for each $p$ the
restriction $K^{p,*}\nach K(X)^{p,*}$ is a homotopy equivalence
over $R^0$.
So we see that there is a well defined map of DG-algebras
$\pi\ot 1:K\ot_RA\nach K(X)\ot_{\tilde{R}}A$, which is by
the properties (c) and (d) even a homomorphism of DDG-algebras.
Further we see that for each $p$ the restriction
$(K\ot_RA)^{p,*}\nach (K(X)\ot_{\tilde{R}}A)^{p,*}$
is a homotopy equivalence.
Hence $\pi\ot 1$ induces a quasi-isomorphism
\begin{align*}
S\ot_RA=\tot(K)\ot_RA=\tot(K\ot_RA)\nach\\
\tot(K(X)\ot_{\tilde{R}}A)=\tot(K(X))\ot_{\tilde{R}}A=
\tilde{S}\ot_{\tilde{R}}A.
\end{align*}
But a quasi-isomorphism of free algebras is already
a homotopy-equivalence.
\qed

Again, we have to explain that the same construction also works
in the simplicial case:
So suppose that $k$ is an object of $\grcn$ and $A$ is
a free algebra over $k$ in $\grcn$.
Say $A=k\lan T\ran$, where each $t\in T$ is associated to a pair
$(\alpha_t,\tau_t,g_t)\in \sN\times\T\times\ZZ_{\leq 0}$.
Now $T^+$ and $X:=T^-$ are sets of free generators in the simplicial sense.
Write $X=\{x_i:\;i\in I\}$ and $(\alpha_i,\tau_i,g_i)$ for the triple
associated to $x_i$.
For $\alpha\in\sN$ set
$$X_\alpha:=\{\rho_{\alpha\alpha_i}(x_i):\;\alpha_i\sub\alpha\}.$$
In the sequel we will simply write $x_i$ for the
element $\rho_{\alpha\alpha_i}(x_i)$ of $R_\alpha$.
Let $E=\{e_i:\;i\in I\}$ be a family of free algebra generators
containing for each $x_i\in X$ an $e_i$ of degree $(g(x),-1)$, belonging
to the simplex $\alpha_i$.
We form the free algebra $K(X)=R\lan E\ran$ in $\gr^2(\sC)^{\sN}$.
Set $E_\alpha:\{\rho_{\alpha\alpha_i}(e_i):\;\alpha_i\sub\alpha\}$.
Then we have $K(X)_\alpha=R_\alpha\lan E_\alpha\ran$.
For each $x=x_i$ in $X$ we fix a family $a_{xy},y\in Y$ with
$Y\sub X_{\alpha_i}$ and $a_{xy}\in R_{\alpha_i}$, such that
$s_\alpha(x)=\sum_ya_{xy}y$.

\begin{prop}
There is a g-finite family $F=\{f_j:\;j\in J\}$ of free algebra generators,
where $f_j$ belongs to $\alpha_j$ and is of bidegree $(g_j,-1)$
and a DDG-algebra structure $(L,h,v)$ on
$L=R\lan E\cup F\ran$ over $R$, such that for all $\alpha$
and all $x\in X_\alpha$ and all
$f\in F_\alpha:=\{f_j:\;\alpha_j\sub\alpha\}$  the
following conditions hold:
\begin{liste}
\item\label{emmq}
$v_{\alpha}(e(x))=x$
\item
$v_\alpha(f)=0$
\item
$h_\alpha(f)\in(X_\alpha\cup F_\alpha)$
\item\label{emmr}
$h_\alpha(e(x))=\sum_ya_{xy}e(y)+\gamma(x)$ for a $\gamma(x)$ in
$(F_{\alpha})$.
\end{liste}
\end{prop}
\bew
We reduce the proposition by induction on the following statement:
Suppose that there is a family $F^{(n)}$ and a DDG-algebra structure
on $(R\lan E\ran_\alpha)_{|\alpha|\leq n}\lan F^{(n)}\ran$, such that
the conditions (i)-(iv) hold for all $\alpha\in\sN^{(n)}$.
Then there is a family
$F^{(n+1)}$ and a DDG-algebra structure on
$(R\lan E\ran_\alpha)_{|\alpha|\leq n+1}\lan F^{(n+1)}\ran$
such that the conditions (i)-(iv) 
hold for all $\alpha\in\sN^{(n+1)}$.
The case $n=0$ as well as the induction step can be done easily
as in the affine case.
\qed

\begin{prop}
There is a g-finite family $G=\{g_j:\;j\in J\}$ o
of free algebra generators,
where $g_j$ belongs to $\alpha_j$ and is of bidegree $g_j$
and a DDG-algebra structure on $K=L\lan G\ran$ over $L$, such that for all
$\alpha\in\sN$ and all $x\in X_\alpha$ and all $g\in G_\alpha$
the following conditions hold:
\begin{liste}
\item
$v_\alpha(g)\in(G_\alpha\cup F_\alpha)$
\item
$h_\alpha(g)\in (G_\alpha).$
\end{liste}
\end{prop}
\bew
With the same method as above, we reduce the statement to the affine
case.
\qed

Now we see that theorem~\ref{Adt}.\ref{hohoho} 
holds as well in the simplicial context.\\

\subsection{A HKR-type theorem}
When we now resume what we know about the resolution
$K(X)\ot_{\tilde{R}}A$, we get the following result.
It generalizes in a sense
the classical Hochschild-Kostant-Rosenberg theorem. 
We use
the same assumptions as in the last subsection.

\begin{satz}\label{HKR-typ}
Consider a homomorphism
$k\nach a$ in $\sC^{\sN}$. Suppose that $\QQ\sub a$.
Let $A$ be a resolvent of $a$
in $\grc^{\sN}(k)$.
Then there is a quasi-isomorphism
\begin{equation*} 
\tot(\wedge_{\cdot}\Omega_{A/k})\nach\Hk(a/k)
\end{equation*}
of DG-algebras in $\grc^{\sN}$ over $A$ and a quasi-isomorphism
\begin{equation*}
\tot(\dach\LL_{a/k})\nach\Hk(a/k)
\end{equation*}
in $\grcn$ over $a$. 
\end{satz}
\vspace{-0.4cm}
\begin{kor}
When $a$ is already free over $k$ (in this case there is no
need to assume that $\QQ\sub a$) and $A=a$, then
$\Omega_{a/k}$ is an object of $\sC^{\sN}_a$ and we get isomorphisms
\begin{equation*}
\wedge^n_a\Omega_{a/k}\isom H_n(\Hk(a/k))
\end{equation*}
Dually, with $T_{A/k}:=\Hom_A(\Omega_{A/k},A)$ we get
$$H^n(\Hom_a(\Hk(a/k),a))\isom\dach^nT_{a/k}.$$
\end{kor}
\bew
The first statement follows directly from the theorem. For the second
statement, remark that in the case, where $A=a$, we
have a quasi-isomorphism of free DG-algebras in $\grcn$ over $A$:
$$\wedge^n_a\Omega_{a/k}\nach\Hk(a/k),$$
where the differential on the left side is trivial.
With proposition~\ref{Apoc}.\ref{rehoeq} it is even a homotopy equivalence.
So with \cite{BinKos} lemma (3.7), the dual homomorphism 
\begin{align*}
\Hom_a(\Hk(a/k),a)\nach\Hom_a(\dach_a\Omega_{a/k},a)=\dach^\cdot T_{a/k} 
\end{align*}
is also a homotopy equivalence, which proves the statement.
\qed

\subsection{The decomposition theorem}

\begin{satz}
We have the following decomposition of Hochschild cohomology:
\begin{equation*}
\HH^n(a/k,M)=\coprod_{i-j=n}H^i(\Hom_A(\dach_A^j\Omega_{A/k},M))
\end{equation*}
\end{satz}
\bew
\begin{align*}
\HH^n(a/k,M)=H^n(\Hom_a(\Hk(a/k),M))=&H^n(\Hom_A(S\ot_RA,M))=\\
H^n(\Hom_A(\tot(K(X)\ot_{\tilde{R}} A),M))=&
H^n(\Hom_A(\tot(\dach^{\cdot}\Omega_{A/k}),M))=\\
H^n(\Hom_A(\coprod_{j\geq 0}\dach^j_A\Omega_{A/k}[j],M))=&
H^n(\prod_j\Hom_A(\dach_A^j\Omega_{A/k}[j],M))=\\
\prod_{j\geq 0}H^{n+j}(\Hom_A(\wedge^j\Omega_A,M))=&
\coprod_{i-j=n}H^{i}(\Hom_A(\dach_A^j\Omega_{A/k},M)).
\end{align*}
The first equality holds by definition. The second one follows by
equation~\ref{glei}. The third one follows by
the simplicial version of theorem~\ref{Adt}.\ref{hohoho}.
The fourth equality holds by remark~\ref{Adt}.\ref{til}.
The fifth follows by equation~\ref{totto}.
The other equalities are elementary.
\qed

\section{Application to complex spaces and noetherian schemes}\label{Appl}
\newcommand{\ch}{\check}
In this section, all schemes and complex spaces are supposed to
be paracompact and separated. For details on many of the 
constructions we refer to \cite{BuFl1} and \cite{BuFl2}.\\

First we will sketch the correlation between the
the theory of coherent sheaves on schemes or complex spaces and
the theory of $\sN$-objects in good pairs of categories.
The main tools that we need here are:
\begin{enumerate}
\item
Instead of considering a space $X$, we consider the simplicial scheme,
associated to an affine covering of $X$. By an affine subspace, we mean
an open affine subscheme in the case of schemes and a 
Stein compact in the case of complex spaces.
There are functors that make simplicial modules out of
sheaves of modules and functors that do the inverse.
\item\label{ekat}
For  affine subsets $U\sub X$ we use the equivalence of categories,
of coherent $\Oh_U$-modules and finite modules
over the ring $\Gamma(U,\Oh_X)$. (Remember that $\Gamma(U,\Oh_X)$
is noetherian, when $X$ is an analytic space.) This equivalence is
given by Cartans theorem A in the analytic case and by
\cite{Hart} chapter 2, exc. 2.4 in the algebraic case.
\end{enumerate}

Now, more generally, let $X$ be a ringed space and $(X_i)_{i\in I}$
a covering of $X$. The \textbf{Nerf} $\sN$ of this covering is the
set of all subsets $\alpha\sub I$, such that 
$\cap_{i\in \alpha}X_i\neq\emptyset$. $\sN$ 
is a simplicial scheme in the sense of subsection~\ref{sims}. 
Further there is a contravariant functor from $\sN$ in the
category of ringed spaces, mapping an object $\alpha$
to the object $\X_\alpha:=\cap_{i\in\alpha}X_i$. 
For $\alpha\sub\beta$ denote the inclusion $X_\beta\nach X_\alpha$
by $p_{\alpha\beta}$.
Such a functor
is called \textbf{simplicial scheme of ringed spaces}. 
Let $X_\ast=(X_\alpha)_{\alpha\in \sN}$ be a simplicial scheme
of ringed spaces. 
Now we define the category of $\Oh_{X_\ast}$-modules.
Its objects
are the families $\sF_\ast=(\sF_\alpha)_{\alpha\in \sN}$ 
with $\sF_\alpha$ in $\MOD(X_\alpha)$ together
with compatible maps $p_{\alpha\beta}^\ast\sF_\alpha\nach\sF_\beta$.
For $\Oh_{X_\ast}$-modules $\sF,\sG$, we set
$\Hom_{X_\ast}(\sF,\sG)$ to be the set of families 
$f_\alpha:\sF_\alpha\nach\sG_\alpha$, that are compatible.
We denote this category by $\MOD(X_\ast)$. The full subcategory
of those $\sF_\ast$, where each $\sF_\alpha$ is coherent is 
denoted by $\COH(X_\ast)$.

\begin{defi} Let $A$ and $B$ be simplicial schemes over
the index sets $A_0$ and $B_0$.
Suppose that $X_\ast=(X_\alpha)_{\alpha\in A}$ and 
$Y_\ast=(Y_\beta)_{\beta\in B}$ are simplicial schemes of ringed spaces.
A \textbf{morphism} $f:X_\ast\nach Y_\ast$ consists of
a mapping $\tau:A_0\nach B_0$, such that for $\alpha\in A$ we get
$\tau(\alpha)\in B$, and a family of compatible maps
$f_\alpha:X_\alpha\nach Y_{\tau(\alpha)}$.
\end{defi}

As in \cite{Flen}, we can form the adjoint functors
\begin{align*}
f^\ast:&\MOD(Y_\ast)\nach\MOD(X_\ast)\quad\text{ and}\\
f_\ast:&\MOD(X_\ast)\nach\MOD(Y_\ast).
\end{align*}
For $\sF$ in $\MOD(Y_\ast)$ and $\alpha\in A$, we have
$(f^\ast\sF)_\alpha:=f_\alpha^\ast\sF_{\tau(\alpha)}$.
The construction of $f_\ast$ is more complicated. For the general
case, we refer to \textit{loc.cit.} But we need only the
following special case:

\begin{bem}
Let $\sF_\ast$ be an object of $\MOD(X_\ast)$. Then for elements
$\beta\in B$ of the form $\beta=\tau(\alpha)$, we have
$$(f_\ast\sF)_\beta=f_{\alpha\ast}\sF_\alpha.$$
\end{bem}
Hence, if the map $\tau:A_0\nach B_0$ is surjective, then the 
construction of $f_\ast$ becomes very simple.

\begin{beisp}
\begin{liste}
\item
When $X$ is a scheme or a complex space and $(X_i)_{i\in I}$
is a covering by affine subspaces, 
then by the separated condition, all $X_\alpha$ are affine.
Now let $(\sC,\sM)$ be the good pair
$(\sC^{(0)},\sM^{(0)})$ or $(\sC^{(1)},\sM^{(1)})$ (see~\ref{Apoc}.\ref{exop}).
Then $a_\ast:=(\Gamma(X_\alpha,\Oh_{X_\alpha}))_{\alpha\in \sN}$ is an
$\sN$-object in $\sC$ and 
there is a 1:1-correspondence between
the objects of $\COH(X_\ast)$ and the $\sN$-objects $M_\ast$ in $\sM$
over $a_\ast$, such that each $M_\alpha$ is finite over $a_\alpha$.
\item
When $X$ is a complex space, and the covering $(X_i)_{i\in I}$ 
is locally finite and chosen in such 
a way that each $X_i$ admits a closed embedding into a polydisc
$P_\alpha$, then we get another simplicial scheme of Stein compacts:
Set $P_\alpha:=\prod_{i\in\alpha} P_i$. Then for $\alpha\sub\beta$,
we have the projection $P_\beta\nach P_\alpha$.
This makes $P_\ast=(P_\alpha)_{\alpha\in\sN}$ a simplicial scheme of
Stein compacts and there is a closed embedding $X_\ast\nach P_\ast$.
\item
Let $X$ be a scheme of finite type over a Ring $\KK$ and $(X_i)_{i\in I}$
an open affine covering of $X$. Again, we can construct a new simplicial 
scheme: Set $a_\alpha:=\Gamma(X_\alpha,\Oh_{\X_\alpha})$ for
$\alpha\in \sN$.
For each $\alpha$, there is a free, finitely generated algebra
$\KK[T]$, that maps surjectively onto $a_\alpha$. So we get a closed
embedding $X_\alpha\nach\spec(\KK[X])=:P_\alpha$.
As above, we get a simplicial scheme $P_\ast$ and a closed embedding
$X_\ast\nach P_\ast$.
\end{liste}
\end{beisp}
The inclusions $j_\alpha:X_\alpha\nach X$ give rise to a map
$j:X_\ast\nach X$ of simplicial schemes of ringed spaces.
Next we will study the adjoint functors $j_\ast$ and $j^\ast$:\\
$j^\ast$ is just the exact functor, mapping an $\Oh_X$-module $\sF$
to the $\Oh_{X_\ast}$-module $(\sF|_{X_\alpha})_{\alpha\in \sN}$.
To describe $j_\ast$, we consider the Cech-functor:
For an $\Oh_{X_\ast}$-module $\sM_\ast$ set
$$\ch{C}^p(\sF):=\prod_{|\alpha|=p}j_{\alpha\ast}\sF_\alpha$$
and define a differential on $\ch{C}^\bullet(\sF_\ast)$ in the usual sense.
Then $j_\ast\sF_\ast$ is just
$H^0(\ch{C}^\bullet(\sF_\ast))$.\\

$j_\ast j^\ast$ is the identity functor.
One can prove the adjointness of $j^\ast$ and $j_\ast$ directly
by a gluing argument. Since $j^\ast$ is an exact functor and
$j_\ast$ is right adjoint to $j^\ast$, we see that $j_\ast$ 
transforms injective objects in $\MOD(X_\ast)$
into injective objects in $\MOD(X)$. \\

By \cite{BuFl1}, proposition 2.26, each $\Oh_{X_\ast}$-module admits 
an injective resolution by modules of the form
$\prod_{\alpha\in \sN}p_{\alpha\ast}\sI_\alpha$ with injective
$\Oh_{X_\alpha}$-modules $\sI_\alpha$.
We will use the following properties of the functor $\ch{C}^\bullet$:
\begin{bem}
\begin{liste}
\item
For $p\geq 0$, the functor $\ch{C}^p$ is exact.
\item
If  $\sF_\alpha$ is an $\Oh_{X_\alpha}$-module, then 
$\ch{C}^\bullet(p_{\alpha\ast}\sF_\alpha)$ is a resolution of 
$j_\ast(p_{\alpha\ast}\sF_\alpha)$.
\item
If $\sF$ is an $\Oh_X$-module, then $\ch{C}^\bullet(j^\ast\sF))$ is a 
resolution of $\sF$.
\end{liste}
\end{bem}

We generalize a part of \cite{BuFl1}, proposition 2.28 for the case where
$X$ is just a ringed space and $X_\ast$ is the simplicial scheme
of ringed spaces associated to some covering $(X_i)_{i\in I}$ of $X$:

\begin{prop}
The functor $j^\ast:D(X)\nach D(X_\ast)$ embeds $D(X)$ as a full and exact
subcategory into $D(X_\ast)$ and $\ch{C}^\bullet=Rj_\ast$ is an exact
right adjoint.
In particular, for $\sF,\sG\in D(X)$ and $\sM_\ast\in D(X_\ast)$, there
are functorial isomorphisms
\begin{align*}
\Ext_X^k(\sF,\sG)&\isom\Ext_{X_\ast}^k(j^\ast\sF,j^\ast\sG)\quad\text{and}\\
\Ext_{X_\ast}^k(j^\ast\sF,\sM_\ast)&\isom\Ext_X^k(\sF,\ch{C}^\bullet(\sM)).
\end{align*}
When all the maps $p_{\alpha\beta}^\ast(\sM_\alpha)\nach\sM_\beta$
for $\alpha\sub\beta$ in $\sN$ are quasi-isomorphisms, then for all n,
there are isomorphisms
\begin{equation*}\label{99021}
\Ext^n_{X_\ast}(\sM_\ast,j^\ast\sF)\isom\Ext^n_X(\ch{C}^\bullet(\sM_\ast),\sF).
\end{equation*}
\end{prop}
\bew
For the proof, that $\ch{C}^\bullet$ is the right derived functor of $j_\ast$, 
we use an injective resolution $\sI_\ast$
of an $\Oh_{X_\ast}$-module $\sF_\ast$ from the same form as above.
Then we have
\begin{align*}
(Rj_\ast)(\sF_\ast)=(j_\ast\sI_\ast)^\bullet=
\prod j_\ast(p_{\alpha\ast}\sI_\alpha)^\bullet\approx
\prod \ch{C}^\bullet(p_{\alpha\ast}\sI_\alpha)=
\ch{C}^\bullet(\sI_\ast^\bullet)=\ch{C}^\bullet(\sF_\ast).
\end{align*}
We only prove the first formula for $\Ext$. Here $\sI_\ast^\bullet$
denotes an injective resolution of $j^\ast\sG$.
\begin{align*}
\Ext^n_{X_\ast}(j^\ast\sF,j^\ast\sG)=
H^n(\Hom_{X_\ast}(j^\ast\sF,\sI_\ast^\bullet))=
H^n(\Hom_X(\sF,j_\ast\sI_\ast^\bullet))=\quad\quad\\
\Ext^n_X(\sF,j_\ast\sI_\ast^\bullet)=\Ext^n_X(\sF,(Rj_\ast)(j^\ast\sG))=\quad\quad\\
Ext_X^n(\sF,\ch{C}^\bullet(j^\ast\sG))=\Ext^n_X(\sF,\sG).\quad\quad\square
\end{align*}

In the sequel, let $\X$ be a complex space or a
scheme of finite type over a noetherian ring.

Now the structure sheaf $\Oh_{\X}$ defines an $\sN$-Object 
$a=a_{\ast}$ in $\sC$.
In the algebraic case each $\Oh_{\X}$-module $\sF$ defines an $\sN$-object
$F=F_{\ast}$ in $\sM$ over $a$.
In the analytic case each coherent $\Oh_{\X}$-module $\sF$ 
defines an $\sN$-object
$F=F_{\ast}$ in $\sM$ over $a$.
Here $(\sC,\sM)$ stands for $(\sC^{(0)},\sM^{(0)})$ in the algebraic case 
and for $(\sC^{(1)},\sM^{(1)})$ in the analytic case.\\

We make the following convention to avoid the distinction between
analytic and algebraic tensor products:\\

\textbf{Convention:} Let $f:X_\ast\nach Y_\ast$ be a morphism of
simplicial schemes of Stein compacts
and let $\sF,\sG$ be graded objects in $\MOD(X_\ast)$, coherent 
in each degree. Then by $\sF\ot_{\Oh_Y}\sG$, we mean the object
in $\MOD(X_\ast)$, which is given by the sheafification
of the object $T_\ast$ in $\grcn$ given as follows:\\
For $\alpha\in\sN$ set 
$B_\alpha:=\Gamma(Y_{\tau(\alpha)},\Oh_{Y_{\tau(\alpha)}})$,
$F_\alpha:=\Gamma(X_\alpha,\sF_\alpha)$ and
$G_\alpha:=\Gamma(X_\alpha,\sG_\alpha)$.
Then $F_\alpha$ and $G_\alpha$ are modules over $B_\alpha$ via the comorphism
of $f_\alpha$. Set $T_\alpha:= F_\alpha\ot_{B_\alpha}G_\alpha$
This defines a simplicial DG algebra $T_\ast$.\\

In the same manner, we define the tensor product $\sF\ot_{\sR}\sG$, when
$\sF$ and $\sG$ are modules over a sheaf of $\Oh_{\X_\ast}$-modules
$\sR$, coherent in each degree.

\subsection{Hochschild-cohomology for complex spaces and schemes}\label{hacop}

Let $f:X\nach Y$ be a morphism of complex spaces
or a morphism of finite type of noetherian schemes.\\

By a \textbf{resolvent} of $X$ over $Y$, we understand a collection of the 
following things:\\
(1) The simplicial scheme $Y_\ast$ associated to a local finite affine covering
$(Y_j)_{j\in J}$ of $Y$; 
(2) the simplicial scheme $X_\ast=(X_\alpha)_{\alpha\in\sN}$ associated to
a local finite affine covering $(X_{ji})_{j\in J,i\in I_j}$
of $X$. This covering is chosen in a way such that for a fixed $j\in J$
the family $(X_{ji})_{i\in I_j}$ is a covering of $f^{-1}(Y_j)$;
(3) a simplicial scheme $P_\ast=(P_\alpha)_{\alpha\in\sN}$ with the same index 
category;
(4) a commutative diagram 
\begin{equation*}
\xymatrix{
X_\ast\ar[r]^{\iota}\ar[d]^{\bar{f}}& P_\ast\ar[dl]^g \\
Y_\ast
}
\end{equation*}
Here $\bar{f}=(\bar{f},\tau)$ is the induced map of simplicial schemes,
$\iota$ is a closed embedding and $g$ is a smooth map\footnote{
This means that for each $\alpha\in\sN$ and each $p\in P_{\alpha}$
the stalk $\Oh_{P_\alpha,p}$ is free (in the analytic case as local
analytic algebra) over $\Oh_{Y_{\tau(\alpha)},y}$.};
(5) a free resolution $\sA_\ast$ of $\Oh_{X_\ast}$ as sheaf of DG-algebras
on $P_\ast$ with $\sA^0_\ast=\Oh_{P_\ast}$, such in each degree there are
only a finite number of free algebra generators.\\

If $\sA_\ast\nach\sB_\ast$ is a morphism of sheaves of DG-algebras,
coherent in each degree, on
a simplicial space $X_\ast$, where each $X_\alpha$ is affine,
then going over to global sections, we can construct a free resolution $S_\ast$
of $B_\ast:=(\Gamma(X_\alpha,\sB_\alpha))_{\alpha\in\sN}$ over
$A_\ast:=(\Gamma(X_\alpha,\sA_\alpha))_{\alpha\in\sN}$,
at least when $B^0_\ast$ is a quotient of a free algebra over $A_\ast^0$
in $\grc^{\sN}$. This follows by
\cite{BinKos}, prop. 8.8. Sheafifying $S_\ast$, we get a free resolution $\sS_\ast$ of
$\sB$ over $\sA$. 
Using this remark, it is easy to deduce the existence of 
free resolutions in the cases we are going to consider.
\begin{beisp}\label{mfkten}
Consider the case where $X$ is smooth and $Y$ is just the single
point $\spec(\CC)$.
Here we can choose $P_i=X_i$ for $i$ in the index set $I$.
Then $X_\alpha$ is a diagonal in $P_\alpha$ and
$A$ can be chosen
to be a Koszul resolution of 
$a=(\Gamma(X_\alpha,\Oh_{X_\alpha}))_{\alpha\in\sN}$ over
$A^0=(\Gamma(P_\alpha,\Oh_{P_\alpha}))_{\alpha\in\sN}$.
In this case one can prove that for each $\alpha$,
$\Omega_{A_\alpha}$ is a module
resolution of $\Omega_{a_\alpha}$.
It follows, that for $\alpha\sub\beta$,
the restriction maps $\LL_\alpha(a/\CC)\nach\LL_\beta(a/\CC)$
are quasi-isomorphisms.
Consequently, the canonical map
$\LL(X)\nach\Omega_X$ is a quasi-isomorphism. 
\end{beisp}

Now, suppose that there is a given resolution
$(X_\ast,Y_\ast,P_\ast,\sA_\ast)$ of the morphism
$f:X\nach Y$. Furthermore, set $\sR:=\sA\ot_{\Oh_{Y_\ast}}\sA$ and
let $\sS$ be a free resolution of $\sA$ over $\sR$.

The following definition coincides for complex spaces with the one,
given in \cite{BuFl2}:

\begin{defi}\label{Hok}
The \textbf{simplicial Hochschild complex} of $X$ over $Y$
is the object in the derived category $D(X_\ast)$ of $\Oh_{X_\ast}$-modules,
represented by 
$$\Hk_\ast(X/Y):=\sS\ot_{\sR}\Oh_{X_\ast}.$$
The \textbf{Hochschild complex} of $X$ over $Y$ is defined as
the object in $D(X)$, represented by
$$\Hk(X/Y):=\ch{C}^\bullet(\Hk_\ast(X/Y)).$$
When $Y$ is just the simple point, we will write
$\Hk(X)$ instead of $\Hk(X/Y)$.
\end{defi}

To show the independence of the Hochschild complex of the choice of the
resolvent, we have to use the following version of \cite{BinKos}, lemma (13.7):
\begin{lemma}\label{lebnd}
Let $f:X\nach X'$ be a flat homomorphism of complex spaces resp.
schemes and
$(X_i)_{i\in I}$ and $(X'_i)_{i\in I'}$ be compact 
locally finite coverings of $X$ and $X'$
by Stein compacts resp. open affine subsets. Let $\tau:I\nach I'$ be a mapping,
such that $f(X_i)\sub X_{\tau(i)}$ for all $i\in I$. Denote the associated
simplicial schemes by $X_\ast$ and $X'_\ast$. Then $f$ defines
a homomorphism $(\bar{f},\tau)$ of simplicial schemes of ringed spaces.
Let $\sG^\bullet$ be a complex in $\COH(X')$ such that for
$\alpha\sub\beta$ the restriction map
$p_{\alpha\beta}^\ast\sG^\bullet_\alpha\nach\sG^\bullet_\beta$
is a quasi-isomorphism.
Then the canonical homomorphism
$$f^\ast\ch{C}(\sG^\bullet)\nach\ch{C}(\bar{f}^\ast\sG^\bullet)$$
is a quasi-isomorphism.
\end{lemma}
\begin{prop}
The definition of $\Hk(X/Y)$ depends 
neither on the resolvent $(Y_\ast,X_\ast,P_\ast,\sA_\ast)$
nor on the choice of
the resolvent $\sS$. 
\end{prop}

\bew
Let $(Y_\ast,X_\ast,P_\ast,\sA_\ast)$ and
$(\tilde{Y}_\ast,\tilde{X}_\ast,\tilde{P}_\ast,\tilde{\sA}_\ast)$
be two resolvents, $\sS$ a free resolution of $\sA$ over $\sA\ot\sA$
and $\tilde{\sS}$ a resolvent of $\tilde{\sA}$ over
$\tilde{\sA}\ot\tilde{\sA}$.
We have to show that there is a quasi-isomorphism
$$\ch{C}(\tilde{\sS}\ot_{\tilde{\sR}}\Oh_{\tilde{X}_\ast})\nach
\ch{C}(\sS\ot_{\sR}\Oh_{X_\ast}).$$
\textbf{First case:} Suppose that $Y_\ast=\tilde{Y}_\ast$,
$X_\ast=\tilde{X}_\ast$ and $P_\ast=\tilde{P}_\ast$.
Then it follows with proposition~\ref{hchc}.\ref{hhoinv} that there is
a quasi-isomorphism
$$\tilde{\sS}\ot_{\tilde{\sR}}\Oh_{\tilde{X}_\ast}\nach
\sS\ot_{\sR}\Oh_{X_\ast}$$
in $\MOD(X_\ast)$. Applying the Cech functor, this case is proven.\\
\textbf{General case:}
Let $Y'_\ast$ be the simplicial scheme associated to the
covering $\{Y_j\}\cup\{Y'_j\}$ and $X'_\ast$ be the simplicial scheme
associated to the covering $\{X_{ij}\}\cup\{X'_{ij}\}$.
We construct $P'_\ast$ in the canonical way and can find a resolvent
$\sA'$, such that $(Y'_\ast, X'_\ast, P'_\ast, \sA'_\ast)$ forms
another resolvent of $f:X\nach Y$.
There is a commutative diagram
\begin{equation*}
\xymatrix{
X_\ast\ar[r]^h\ar[d]^{f_\ast} & X'_\ast\ar[d]^{f'_\ast}\\
Y_\ast\ar[r] & Y'_\ast }
\end{equation*}
By the first case, there is a quasi-isomorphism
$$h^\ast(\sS'\ot_{\sR'}\Oh_{X'_\ast})\approx\sS\ot_{\sR}\Oh_{X_\ast}.$$
With lemma~\ref{Appl}.\ref{lebnd}, there is a quasi-isomorphism
$$\ch{C}(\sS'\ot_{\sR'}\Oh_{X'_\ast})\approx
\ch{C}(h^\ast(\sS'\ot_{\sR'}\Oh_{X'_\ast})).$$
Hence we get $\ch{C}(\sS\ot_{\sR}\Oh_{X_\ast})\approx
\ch{C}(\sS'\ot_{\sR'}\Oh_{X'_\ast})$ and in the same way we get
$\ch{C}(\tilde{\sS}\ot_{\tilde{\sR}}\Oh_{\tilde{X}_\ast})\approx
\ch{C}(\sS'\ot_{\sR'}\Oh_{X'_\ast})$.
\qed

As in \cite{BuFl2}, we define the \textbf{Hochschild cohomology of
  $\X$ over $\Y$}  
with values in the sheaf $\sF$ 
as $\Ext_{\X_{\ast}}(\Hk(\X/\Y),\sF)$. At least in the case
where $\sF$ is coherent, we want to show that this definition
is equal to the following one, which seems to be more natural,
from the viewpoint of good pairs of categories: 

\begin{defi}\label{99022}[alternative]\\
Suppose that $\sF$ is coherent.
Let $a$ be the algebra $(\Gamma(X_\alpha,\Oh_{X_\alpha}))_{\alpha\in\sN}$
in $\sC^{\sN}$, let $k$ be the algebra
$\Gamma(Y_{\tau(\alpha)},\Oh_{Y_{\tau(\alpha)}}))_{\alpha\in\sN}$ in
$\sC^{\sN}$. Then to $f$, there is associated a
homomorphism $k\nach a$ in $\sC^{\sN}$.
Let $F$ be the module $(\Gamma(X_\alpha,\sF_\alpha))_{\alpha\in\sN}$.
Then we define the n-th \textbf{Hochschild cohomology} of $X$ over $Y$
as
\begin{equation*}
\HH^{n}(\X/\Y,\sF):=H^{n}(\Hom_{a}(\Hk_{\ast}(a/k),F)).
\end{equation*}
\end{defi} 
\begin{bem}\label{hof}
For $\sM_\ast:=\Hk_\ast(\X/\Y)$, the assumption of the second part of 
proposition~\ref{Appl}.\ref{99021} is satisfied, i.e. for
$\alpha\sub\beta$
the maps $p_{\alpha\beta}^\ast(\sM_\alpha)\nach\sM_\beta$ are 
quasi-isomorphisms.
\end{bem}
\bew
\cite{BuFl2}, lemma 1.7.
\qed
\begin{kor}
For coherent $\Oh_{\X}$-modules $\sF$,
the two definitions of Hochschild cohomology coincide, i.e.
\begin{equation*}
\HH^i(\X/\Y,\sF)=\Ext^i_{\X}(\Hk(\X/\Y),\sF).
\end{equation*}
\end{kor}
\bew
Since $\Hk_\ast(\X/\Y)$ is a complex of free $\Oh_{\X_\ast}$-modules,
with proposition~\ref{Appl}.\ref{99021} we get
\begin{align*}
\Ext^i_{\X}(\Hk(\X/\Y),\sF)=\Ext^i_{\X_\ast}(\Hk_\ast(\X/\Y),j^\ast\sF)=\\
H^i(\Hom_a(\Hk_\ast(a/k),F_\ast))=\HH^i(\X/\Y,\sF).
\end{align*}
\qed

For a (noetherian) scheme $X$, we want to show that the definition
of the Hochschild complex $\Hk(X)$ coincides with one of the
definitions given in
\cite{Swan} or \cite{Yeku}:
\begin{prop}
Let $\sC^{\cycl}(X)$ be the complex of sheaves in $\Mod(X)$,
associated to the presheaf $U\mapsto C^{\cycl}(\Gamma(U,\Oh_X))$.
Then in $D(X)$, there is an isomorphism
$$\sC^{\cycl}(X)\isom\Hk(X).$$
\end{prop}
\bew
Choose a resolvent $(X_\ast,P_\ast,\sA_\ast)$ of $X$ over the base
ring $\KK$. Let $\sS$ be a resolvent of $\sA$ over $\sR=\sA\ot\sA$.
Let $a, A, R$ and $S$ be the simplicial algebras in $\grcn$
corresponding
to $\Oh_{X_\ast},\sA,\sR$ and $\sS$.
There is an isomorphism in $D(X)$:
$$\ch{C}(j^\ast\sC^{\cycl}(\Oh_X))\isom C^{\cycl}(\Oh_X).$$
But $j^\ast\sC^{\cycl}(\Oh_X)$ corresponds to $C^{\cycl}(a)$, which  
is quasi-isomorphic to $C^{\ba}(A)\ot_Ra$.
$C^{\ba}(A)$ is a resolution of $A$ over $R$, so by \cite{BinKos},
proposition (8.4), there is a quasi-isomorphism
$f:S\nach C^{\ba}(A)$. Now each $C^{\ba}(A)_\alpha=C^{\ba}(A_\alpha)$
is, as well as $S_\alpha$, a free module resolution of $A_\alpha$ over $R_\alpha$.
So each $f_\alpha$ is a homotopy equivalence.
Hence for each $\alpha$
$$f_\alpha\ot 1:S_\alpha\ot_{R_\alpha}a_\alpha\nach C^{\ba}
(A_\alpha)\ot_{R_\alpha}a_\alpha$$ is a homotopy equivalence.
So $f\ot 1$ is a quasi-isomorphism.
In $D(X_\ast)$ it induces an isomorphism
$$j^\ast(\sC^{\cycl}(\Oh_X))\isom\sS\ot_{\sR}\Oh_X.$$
Forming the Cech complex gives the desired result.
\qed

\subsection{The decomposition Theorem}

The quasi-isomorphism $\tot(\dach\LL_{a/k})\nach\Hk(a/k)$ in $\grm^{\sN}$
over $a$ in theorem~\ref{Adt}.\ref{HKR-typ} defines a quasi-isomorphism
$$\tot(\dach\LL_\ast(X/Y))\nach\Hk_\ast(X/Y)$$ in $\MOD(X_\ast)$.
Since the Cech-functor is exact, we get the following HKR-type
theorem:

\begin{satz}\label{HKRgeo}
There is
an isomorphism
$$\tot(\dach\LL(X/Y))\nach\Hk(X/Y)$$
in the derived category $D(X)$.
\end{satz}

From this we deduce easily the announced decomposition theorem:
 
\begin{kor}
There is a natural decomposition
$$\HH^n(X/Y,\sM)=\oplus_{p-q=n}\Ext_X^p(\tot(\dach^q\LL(X/Y)),\sM)$$
\end{kor}
For complex spaces,
this is just theorem 4.2 of \cite{BuFl2}.
There is another nice description of Hochschild cohomology of
complex spaces or noetherian schemes over a field $K$ in any
characteristic:
\begin{bem}
$HH^n(X)=\Ext_{X^2}(\Oh_X,\Oh_X)$.
\end{bem}
\bew
We will use the letter $\KK$ for the field $K$ or for the complex numbers,
depending on the context. With the notations as above, we get:
\begin{align*}
\HH^n(\X)=H^n(\Hom_A(S\ot_Ra,a))=H^n(\Hom_R(S,a))=\\
H^n(\Hom_{a\ot_{\KK}a}(S\ot_R(a\ot_{\KK}a),a))=
H^n(\Hom_{\Oh_{X^2_\ast}}(\sS\ot_\sR\Oh_{X_\ast^2}))\\ 
=\Ext^n_{\Oh_{\X_\ast^2}}(\Oh_{\X_\ast},\Oh_{\X_\ast})=
\Ext_{X^2}(\Oh_X,\Oh_X).
\end{align*}
Here we have used that $\sS\ot_{\sR}\Oh_{X_\ast^2}$ is a resolution of
$\Oh_{X_\ast}$ and the fact that it is free over $\Oh_{X_\ast^2}$.
\qed

\pagebreak

\subsection{Hochschild-Cohomology for manifolds and smooth 
varieties}   

%\begin{align*}
%H^n(\Hom_A(S\ot_R A,A))=\Hom_A(S\ot_R A,A)=\dach^n\Hom_A(\Omega_{A/k},A)=\\
%=\dach^n T_\X
%\end{align*}

At last we will see that the decomposition theorem for Hochschild
cohomology
of complex manifolds, announced in \cite{Kont}, follows easily from
theorem~\ref{Appl}.\ref{HKRgeo}. It holds as well for smooth
schemes of finite type over a field $\KK$ of characteristic zero.
This case was proven in a different way by Yekutieli \cite{Yeku}.

\begin{satz}
Let $X$ be a complex analytic manifold or a smooth scheme 
of finite type over a field $\KK$ of characteristic zero.
Then there is a decomposition of Hochschild
cohomology: 
\begin{equation*}
\HH^n(X)=\coprod_{i-j=n}H^i(\X,\dach^j\sT_{\X}).
\end{equation*}
\end{satz}
\bew
For complex analytic manifolds,
we work with a fixed covering by Stein compacts
and its associated simplicial scheme
$X_\ast$.
For the case of smooth schemes of finite type over $\KK$,
we work with an open affine covering by schemes of the form
$\spec(A)$, where $A$ is a finitely generated $\KK$-algebra.
Denote the associated simplicial scheme also by $X_\ast$.\\
By proposition~\ref{Appl}.\ref{hof}, theorem~\ref{Adt}.\ref{HKR-typ} 
and example~\ref{Appl}.\ref{mfkten},
there are quasi-isomorphisms
$$j^\ast(\Hk(X))=j^\ast C(\Hk_\ast(X))\approx\Hk_\ast(X)\approx
\dach\LL_\ast(X)\approx
\dach_{\Oh_{X_\ast}}\Omega_{X_\ast}=j^\ast(\dach_{\Oh_X}\Omega_X)$$
of $\Oh_X$-modules.
$j_\ast j^\ast$ is the identity functor, so there is a 
quasi-isomorphism
of $\Oh_X$-modules
$$\Hk(X)\approx\dach_{\Oh_X}\Omega_X.$$

We consider $\wedge_{\Oh_X}\Omega_X$ as complex in negative degrees,
so $\wedge\Omega_X=\coprod_{j\geq 0}\wedge^j\Omega_X[j]$
and
$$\HH^n(X)=\Ext^n_X(\Hk(X),\Oh_X)\isom
\coprod_{j\geq 0}\Ext_X^{n+j}(\wedge^j\Omega_X,\Oh_X).$$
By \cite{Gode}, theorem (7.3.3),  there is a (bounded)
spectral sequence with
$E_2^{p,q}=H^p(X,\EXT_X^q(\wedge^j\Omega_X,\Oh_X))$, converging to
$\Ext_X(\wedge^j\Omega_X,\Oh_X)$.
But $\dach^j\Omega_X$ is a locally free $\Oh_X$-module, so 
$\EXT_X^q(\dach^j\Omega_X,\Oh_X)$ is zero for $q>0$
and $\HOM_\X(\dach\Omega_X,\Oh_X)$ for $q=0$.
So the spectral sequence degenerates at once
and we get
$$\Ext_X^q(\wedge^j\Omega_X,\Oh_X)=
H^q(X,\HOM_X(\wedge^j\Omega_X,\Oh_X).$$
Now 
there is a natural isomorphism of sheaves
$$\dach^j\sT_X=\dach^j\HOM_X(\Omega_X,\Oh_X)
\nach\HOM_X(\dach^j\Omega_X,\Oh_X),$$
which by \cite{BAlg3}, prop. 7, p. 154 is an isomorphism. 
So we get
$\HH^n(X)=\coprod_{j\geq 0}H^{n+j}(X,\wedge^j\sT_X)=
\coprod_{i-j=n}H^i(X,\wedge^j\sT_X)$.
\qed

\lz
\begin{center}
Institut Fourier\\
UMR 5582\\
BP 74\\
38402 Saint Martin d'H\`eres\\
France\\
\lz
frank.schuhmacher@ujf-grenoble.fr
\end{center}

\end{document}